\documentclass[12pt,a4paper,reqno]{article}%
\usepackage{amssymb}
\usepackage{amsmath}
\usepackage{amsfonts}
\usepackage{graphicx}
\usepackage{bm}
\usepackage{mathrsfs}
\usepackage[all]{xy}
\usepackage{titletoc}
\usepackage{rotating}
\usepackage[T1]{fontenc}%
\SelectTips{cm}{}
\newtheorem{theorem}{Theorem}
\newtheorem{corollary}[theorem]{Corollary}
\newtheorem{definition}[theorem]{Definition}
\newtheorem{lemma}[theorem]{Lemma}
\newtheorem{proposition}[theorem]{Proposition}
\newtheorem{remark}[theorem]{Remark}
\newtheorem{example}[theorem]{Example}
\newenvironment{proof}[1][Proof]{\noindent\textbf{#1.} }{\ \rule{0.5em}{0.5em}}
\title{\bf Secondary Calculus and\\ the Covariant Phase Space}
\author{\sc{L.~Vitagliano}\thanks{{\bf e}-{\it mail}: \texttt{lvitagliano@unisa.it}}\\
\small{DMI, Universit\`a degli Studi di Salerno, and}\\ \small{Istituto Nazionale di Fisica Nucleare, GC Salerno}\\
\small{Via Ponte don Melillo, 84084 Fisciano (SA), Italy}}

\DeclareFontFamily{OT1}{wncyi}{}
\DeclareFontShape{OT1}{wncyi}{m}{it}{
<5> <6> <7> <8> <9> gen * wncyi
<10> <10.95> <12> <14.4> <17.28> <20.74> <24.88> wncyi10
}{}
\DeclareSymbolFont{cyrletters}{OT1}{wncyi}{m}{it}
\DeclareSymbolFontAlphabet{\cyrmath}{cyrletters}
\DeclareMathSymbol{\rE}{\cyrmath}{cyrletters}{003}
\DeclareMathSymbol{\rD}{\cyrmath}{cyrletters}{068}
\DeclareMathSymbol{\rG}{\cyrmath}{cyrletters}{017}
\DeclareMathSymbol{\rI}{\cyrmath}{cyrletters}{073}
\DeclareMathSymbol{\rL}{\cyrmath}{cyrletters}{076}
\DeclareMathSymbol{\rZ}{\cyrmath}{cyrletters}{090}

\begin{document}
\maketitle
\begin{abstract}
The covariant phase space of a Lagrangian field theory is the solution space of the associated Euler-Lagrange equations. It is, in principle, a nice environment for covariant quantization of a Lagrangian field theory. Indeed, it is manifestly covariant and possesses a canonical (functional) \textquotedblleft presymplectic structure\textquotedblright\ $\boldsymbol{\omega}$ (as first noticed by Zuckerman in 1986) whose degeneracy (functional) distribution is naturally interpreted as the Lie algebra of gauge transformations. We propose a fully rigorous approach to the covariant phase space in the framework of jet spaces and (A.~M.~Vinogradov's) secondary calculus. In particular, we describe the degeneracy distribution of $\boldsymbol{\omega}$. As a byproduct we rederive the existence of a Lie bracket among gauge invariant functions on the covariant phase space.
\end{abstract}

\newpage
\contentsmargin{2.55em}
\titlecontents{section}
[1.5em]
{}
{\contentslabel{2.3em}}
{\hspace*{-2.3em}}
{\titlerule*[1pc]{.}\contentspage}
\titlecontents{subsection}
[3.8em]
{} % note that 3.8 = 1.5 + 2.3
{\contentslabel{3.2em}}
{\hspace*{-3.2em}}
{\titlerule*[1pc]{.}\contentspage}

\begin{table}[p]
\textbf{List of Main Symbols}\\
\begin{tabular}{ll}
& \\
  $J^{\infty}\pi $ & space of $\infty$--jets of local sections of the bundle
$\pi:E\longrightarrow M$ \\
  $\mathrm{diff}(\pi,\tau) $ & module of differential
operators from $\pi$ to $\tau$ \\

  $\mathscr{E}_{\Phi} $ & PDE determined by the differential operator
$\Phi$ \\

  $\mathscr{E} $ & $\infty$ prolongation of a PDE \\

  $\mathscr{C} $ & Cartan distribution \\

  $\mathscr{C}\mathrm{D}(\mathscr{E}) $ & module of horizontal vector
fields on $\mathscr{E}$ \\

  $\mathscr{C}\Lambda(\mathscr{E}) $ & Cartan ideal of $\mathscr{E}$ \\

  $\mathscr{C}^{p}\Lambda(\mathscr{E}) $ & $p$th exterior power of
$\mathscr{C}\Lambda(\mathscr{E})$ \\

  $\mathscr{C}^{\bullet}\Lambda(\mathscr{E}) $ & algebra generated by
$\mathscr{C}\Lambda^{1}(\mathscr{E})$ \\

  $\mathscr{C}E(\mathscr{E}) $ & $\mathscr{C}$--spectral sequence of
$\mathscr{E}$ \\

  $\overline{\Lambda}(\mathscr{E}) $ & algebra of horizontal forms on
$\mathscr{E}$ \\

  $\overline{d} $ & horizontal de Rham differential \\

  $\overline{H}(\mathscr{E}) $ & horizontal de Rham cohomology of $\mathscr{E}$ \\

  $d^{V} $ & vertical de Rham differential \\

  $V\mathrm{D}(\mathscr{E}) $ & module of
vertical vector fields on $\mathscr{E}$ \\

  $\mathrm{D}_{\mathscr{C}}(\mathscr{E}) $ & Lie algebra of symmetries
of $(\mathscr{E},\mathscr{C})$ \\

  \textrm{Sym}$(\mathscr{E}) $ & Lie algebra of non--trivial symmetries
of $(\mathscr{E},\mathscr{C})$ \\

  $V\mathrm{D}_{\mathscr{C}}(\mathscr{E}) $ & Lie algebra of vertical
symmetries of $(\mathscr{E},\mathscr{C})$ \\

  $\varkappa $ & module of generating sections
of higher symmetries of $\pi$ \\

  $\ell_{\Phi} $ & universal linearization of the differential operator
$\Phi$ \\

  $\boldsymbol{C}^{\infty}(\boldsymbol{M})^{\bullet} $ & space of
secondary functions on the secondary manifold $\boldsymbol{M}$ \\

  $\mathbf{D}(\boldsymbol{M})^{\bullet} $ & space of secondary vector
fields on $\boldsymbol{M}$ \\

  $\boldsymbol{\Lambda}(\boldsymbol{M})^{\bullet} $ & space of secondary
differential forms on $\boldsymbol{M}$ \\

  $\boldsymbol{d} $ & secondary de Rham differential \\

  $\overline{S} $ & horizontal Spencer differential \\

  $\mathscr{C}\mathrm{Diff}(P,Q) $ & module of horizontal differential
operators $P\longrightarrow Q$ \\

  $\overline{J}{}^{\infty}P $ & module of $\infty$ horizontal jets of
elements of $P$ \\

  $\overline{j}{}_{\infty} $ & $\infty$ horizontal jet prolongation
$P\longrightarrow\overline{J}{}^{\infty}P$ \\

  $h_{\square}^{\infty} $ & homomorphism $\overline{J}{}^{\infty
}P\longrightarrow\overline{J}{}^{\infty}Q$ associated to $\square
\in\mathscr{C}\mathrm{Diff}(P,Q)$ \\

  $\eta_{\Phi} $ & natural monomorphism $V\mathrm{D}%
(\mathscr{E})\longrightarrow\overline{J}{}^{\infty}\varkappa|_{\mathscr{E}}$ \\

  $\eta_{\Phi}^{\ast} $ & natural epimorphism $\mathscr{C}\mathrm{Diff}%
(\varkappa|_{\mathscr{E}},\overline{\Lambda}(\mathscr{E}))\longrightarrow
\mathscr{C}\Lambda^{1}(\mathscr{E})\otimes\overline{\Lambda}(\mathscr{E})$ \\

  $\int $ & natural projection $\overline{\Lambda}{}^{n}%
(\mathscr{E})\longrightarrow\overline{H}{}^{n}(\mathscr{E})$ \\

  $\boldsymbol{E}(\mathscr{L}) $ & left hand side of the Euler--Lagrange equations \\

  $\boldsymbol{P} $ & covariant phase space \\

  $\boldsymbol{\omega} $ & canonical, closed, secondary $2$--form on
$\boldsymbol{P}$ \\

  $\Delta_{1} $ & compatibility operator for $\ell_{\boldsymbol{E}%
(\mathscr{L})}$ \\

  $\boldsymbol{\Omega} $ & linear map $\mathbf{D}(\boldsymbol{M}%
)^{\bullet}\longrightarrow\boldsymbol{\Lambda}^{1}(\boldsymbol{M})^{\bullet}$
associated to $\boldsymbol{\omega}$ \\

\end{tabular}

\end{table}

\section*{Introduction}

Covariant phase space (CPS) is the solution space of a system of
Euler--Lagrange partial differential equations\footnote{Notice that sometimes
the name covariant phase space is referred to the quotient of the above
mentioned solution space with respect to gauge transformations.} (PDEs). It
has been first noticed by Zuckerman in the 1986 \cite{z87} (see also
\cite{c88,cw87}) that there is a canonical, closed $2$--form
$\boldsymbol{\omega}$ on such a functional space generalizes the
symplectic form on the phase space of a regular Lagrangian system in
mechanics. Moreover, the degeneracy distribution of $\boldsymbol{\omega}$ is
naturally interpreted as Lie algebra of gauge transformations \cite{lw90}.
Therefore, the CPS is, in principle, a nice environment to perform a covariant
(canonical) quantization of a Lagrangian theory. Namely, gauge invariant
functions on the CPS possess a well defined Lie bracket induced by
$\boldsymbol{\omega}$, which has been proved in \cite{bhs91} to coincide with
the so--called Peierls bracket \cite{p52}. In turn, Peierls bracket is at
the basis of the global approach to quantum field theory \cite{d03}.

Despite its conceptual relevance, the CPS is, in general, a complicated
functional space, which is difficult to handle with analytic methods. Indeed,
most of the literature about it (see \cite{r04} and references therein) comes
from the physicists community and it is rarely completely rigorous from a
mathematical point of view. For instance, it seems to be very hard to
rigorously perform, in full generality, a symplectic reduction of the CPS to
get rid of gauge (non--physical) degrees of freedom.

On the other hand, A.~M.~Vinogradov developed a whole theory, the so--called
\emph{secondary calculus }(see \cite{v01} and references therein, and
\cite{v98} for a short introduction), which properly formalizes in
cohomological terms the idea of a (local) functional differential calculus on
the space of solutions of a generic system of PDEs (for this reason, roughly
speaking, the word \textquotedblleft secondary\textquotedblright\ in this
paper could be considered as a synonym of \textquotedblleft
functional\textquotedblright). Thus, secondary calculus appears to be a
suitable setting to rigorously investigate the CPS and its properties. The aim
of the paper is to describe rigorously the CPS, its canonical $2$--form and
some their properties within secondary calculus. As a byproduct it will become
transparent the analogy between the CPS and the phase space of constrained
mechanical systems.

The paper is divided into two parts. In order to make it as self--consistent
as possible we review, in the first part, those aspects of secondary calculus
that are needed for a suitable formalization of the CPS. In Sections
\ref{Sec1}, \ref{Sec2} and \ref{Sec3} we briefly describe the geometry and the
main properties of jet spaces and differential equations, and relevant
structures on them. In Section \ref{Sec4} we define secondary vector fields
and differential forms, and summarize the main formulas of first order
secondary calculus. In Sections \ref{HorCalc} and \ref{SpencGold} we review
the main technical aspects of secondary calculus and how to handle the
relevant cohomologies.

The second part of the paper is devoted to the CPS and to original results on
the subject. In Section \ref{Sec7} we introduce the CPS for a general
Lagrangian field theory (any number of variable and any order) and rederive
the existence of a canonical $2$--form $\boldsymbol{\omega}$ on it completing
the proof by Zuckerman \cite{z87}. In Section \ref{Sec8} we propose a
\textquotedblleft symplectic version\textquotedblright\ of the first Noether
theorem, which makes it evident the analogy with Hamiltonian mechanics. In
Section \ref{Sec9} we describe the degeneracy distribution of
$\boldsymbol{\omega}$ and propose, and motivate, a new (and very natural)
definition of (infinitesimal) gauge symmetries in field theory. In Section
\ref{Sec10} we describe gauge invariant secondary functions on the CPS and
show that they are endowed with a canonical Lie bracket (such bracket
formalizes rigorously the Peierls bracket \cite{p52}). In Section \ref{Sec11}
we outline a possible path through a \textquotedblleft secondary symplectic
reduction\textquotedblright\ of the CPS. Applications to concrete Lagrangian
theories will be presented somewhere else.

Most of the (almost) trivial computations will be performed in some details to
emphasize similarities between secondary calculus and standard calculus on manifolds.

\subsection*{Notations and Conventions}

In this section we collect notations and conventions about some general
constructions in differential geometry that will be used in the following.

Let $N$ be a smooth manifold. We denote by $C^{\infty}(N)$ the $\mathbb{R}%
$--algebra of smooth, $\mathbb{R}$--valued functions on $N$. We will always
understand a vector field $X$ on $N$ as a derivation $X:C^{\infty
}(N)\longrightarrow C^{\infty}(N)$. The value of $X$ at the point $x\in M$
will be denoted by $X_{x}$. We denote by $\mathrm{D}(N)$ the $C^{\infty}%
(N)$--module of vector fields over $N$, by $\Lambda(M)=\bigoplus_{k}%
\Lambda^{k}(N)$ the graded $\mathbb{R}$--algebra of differential forms over
$N$ and by $d:\Lambda(N)\longrightarrow\Lambda(N)$ the de Rham differential.
If $F:N_{1}\longrightarrow N$ is a smooth map of manifolds, we
denote by $F^{\ast}:\Lambda(N)\longrightarrow\Lambda(N_{1})$ its pull--back.

Let $\alpha:W\longrightarrow N$ be a vector bundle and $F:N_{1}\longrightarrow
N$ a smooth map of manifolds. The $C^{\infty}(N)$--module of smooth sections
of $\alpha$ will be denoted by $\Gamma(\alpha)$. For $s\in\Gamma(\alpha)$ and
$x\in N$ we put, sometimes, $s_{x}:=s(x)$. The zero section of $\alpha$ will
be denoted by $o:N\ni x\longmapsto o_{x}:=0\in\alpha^{-1}(x)\subset W$. The
vector bundle on $N_{1}$ induced by $\alpha$ via $F$ will be denote by
$F^{\circ}(\alpha):F^{\circ}(W)\longrightarrow N$:
\[%
\begin{array}
[c]{c}%
\xymatrix{F^\circ(W) \ar[r] \ar[d]_-{F^\circ(\alpha)} & W \ar[d]^-{\alpha} \\
               N_1 \ar[r]^-F                        &  N }
\end{array}
.
\]
For any section $s\in\Gamma(\alpha)$ there exists a unique section, which we denote by $F^{\circ}(s)\in\Gamma(F^{\circ}(\alpha
))$, such that the diagram
\[
\xymatrix{F^\circ(W) \ar[r]  & W  \\
               N_1 \ar[r]^-F    \ar[u]^-{F^\circ(s)}                    &  N \ar[u]_-{s}}
\]
commutes. If $i_{L}:L\hookrightarrow N$ is the embedding of a submanifold then
we put $\alpha|_{L}:=i_{L}^{\circ}(\alpha)$, $\Gamma(\alpha)|_{L}%
:=\Gamma(\alpha|_{L})$ and for $s\in\Gamma(\alpha)$, $s|_{L}:=i_{L}^{\circ
}(s)$. $s|_{L}$ will be referred to as \emph{the restriction to }$L$ of
$s$.

Let $F:N_{1}\longrightarrow N$ be as above. A vector field along $F$ is an
$\mathbb{R}$--linear map $X:C^{\infty}(N)\longrightarrow C^{\infty}(N_{1})$
such that the following Leibnitz rule holds: $X(fg)=F^{\ast}(f)X(g)+F^{\ast
}(g)X(f)$, $f,g\in C^\infty(N)$. Vector fields along $F$ identify with sections of the induced
bundle $F^{\circ}(\tau_{N}):F^{\circ}(TN)\longrightarrow N_{1}$, $\tau
_{N}:TN\longrightarrow N$ being the tangent bundle to $N$.

Let $\zeta:A\longrightarrow N$ be a fiber bundle. We denote by $\nu
\zeta:V\zeta\longrightarrow A$ the vertical (with respect to $\zeta$) tangent
bundle to $A$ and by $V_{a}\zeta:=(\nu\zeta)^{-1}(a)$ its fiber over $a\in A$.
Notice that $V\zeta\subset TA$, the tangent manifold to $A$. If $\zeta
_{1}:A_{1}\longrightarrow N_{1}$ is another fiber bundle, $F:A_{1}%
\longrightarrow A$ a morphism of fiber bundles and $TF: TA_1 \longrightarrow TA$ the associated tangent map, then $(TF)(V\zeta_{1})\subset
V\zeta$ and, therefore, it is well defined the restriction $VF:V\zeta
_{1}\longrightarrow V\zeta$ of $TF$ to $V\zeta_{1}$ and $V\zeta$, and the
diagram
\[
\xymatrix{V\zeta_1 \ar[r]^{VF} \ar[d]_-{\nu \zeta_1} & V\zeta \ar[d]^-{\nu \zeta} \\
               A_1 \ar[r]^-F                        &  A }
\]
commutes.

Let
\[
\xymatrix{\cdots \ar[r] & K_{l-1} \ar[r]^-{\delta_{l-1}} & K_{l} \ar[r]^-{\delta_{l}} & K_{l+1} \ar[r]^-{\delta_{l+1}} & \cdots}
\]
be a complex. Put $K:=\bigoplus_{l}K_{l}$ and $\delta:=\bigoplus_{l}\delta
_{l}$. We denote by $H(K,\delta):=\bigoplus_{l}H^{l}(K,\delta)$, the
cohomology space of $(K,\delta)$, $H^{l}(K,\delta):=\ker\delta_{l}%
/\operatorname{im}\delta_{l-1}$. If $\omega\in\ker\delta$, then we denote by
$[\omega]$ its cohomology class.

Denote by $\mathbb{N}$ the set of natural numbers and put $\mathbb{N}%
_{0}:=\mathbb{N}\cup\{0\}$. We will always understand the sum over repeated
upper-lower (multi-)indexes. Our notations about multi-indexes are the
following. Let $n\in\mathbb{N}$, $\mathbb{I}_{n}=\{1,\ldots,n\}$ and
$\mathbb{M}_{n}$ be the free abelian monoid generated by $\mathbb{I}_{n}$.
Even if $\mathbb{M}_{n}$ is abelian we keep for it the multiplicative
notation. Thus if $I=i_{1}\cdots i_{l},J=j_{1}\cdots j_{m}\in$ $\mathbb{M}%
_{n}$ are (equivalence classes of) words, $i_{1},\ldots,i_{l},j_{1}%
,\ldots,j_{m}\in\mathbb{I}_{n}$, we denote by $IJ=i_{1}\cdots i_{l}%
j_{1}\cdots j_{m}$ their composition. If $I=i_{1}\cdots i_{l}\in\mathbb{M}%
_{n}$ is a word, $i_{1},\ldots,i_{l}\in\mathbb{I}_{n}$, denote by $|I|:=l$ its
length. We denote by $\mathsf{O}$ the (equivalence class of the) empty word.
An element $I\in\mathbb{M}_{n}$ is called an $n$\emph{-multi-index} (or,
simply, a multi--index) and $|I|$ the length of the multi-index. For $k \leq \infty $ let
$\mathbb{M}_{n,k}\subset\mathbb{M}_{n}$ be the subset made of multi-indexes
of length $\leq k$. If $(x^{1},\ldots,x^{n})$ are local coordinates on a manifold
$N$, $n=\dim N$, and $I=i_{1}\cdots i_{k}\in\mathbb{M}_{n}$, we put
$\tfrac{\partial^{|I|}}{\partial x^{I}}:=\tfrac{\partial^{k}}{\partial
x^{i_{1}}\cdots\partial x^{i_{k}}}$. We stress that this notation is different
from more popular ones (see, for instance, \cite{a92}).

\section{Secondary Calculus}

\subsection{Jet Spaces and PDEs\label{Sec1}}

Let $\pi:E\longrightarrow M$ be a fiber bundle, $\dim M=n$, $\dim E=m+n$. For
$l\leq k\leq\infty$, we denote by $\pi_{k}:J^{k}\pi\longrightarrow M$ the
bundle of $k$-jets of local sections of $\pi$, and by $\pi_{k,l}:J^{k}%
\pi\longrightarrow J^{l}\pi$ the canonical projection. For any local section
$\boldsymbol{p}:U\longrightarrow E$ of $\pi$, $U\subset M$ being an open
subset, we denote by $j_{k}\boldsymbol{p}:U\longrightarrow J^{k}\pi$ its
$k$th jet prolongation and by $\Gamma_{\boldsymbol{p}}^{k}%
:=\operatorname{im}j_{k}\boldsymbol{p}$ its graph. For $x\in U$, put
$[\boldsymbol{p}]_{x}^{k}:=(j_{k}\boldsymbol{p})(x)$. Any system of adapted to
$\pi$ coordinates $(\ldots,x^{i},\ldots,u^{\alpha},\ldots)$ on an open subset
$U$ of $E$ gives rise to a system of jet coordinates $(\ldots,x^{i}%
,\ldots,u_{I}^{\alpha},\ldots)$ on $\pi_{k,0}^{-1}(U)\subset J^{k}\pi$,
$i=1,\ldots,n$, $\alpha=1,\ldots,m$, $I\in\mathbb{M}_{n,k}$, where we put
$u_{\mathsf{O}}^{\alpha}:=u^{\alpha}$, $\alpha=1,\ldots,m$. If a local section
$\boldsymbol{p}$ of $\pi$ is locally given by
\begin{equation}
u^{\alpha}=p^{\alpha}(\ldots,x^{i},\ldots),\quad\alpha=1,\ldots,m, \label{p}%
\end{equation}
then $j_{k}\boldsymbol{p}$ is locally given by
\[
u_{I}^{\alpha}=(\tfrac{\partial^{|I|}}{\partial x^{I}}p^{\alpha})(\ldots
,x^{i},\ldots),\quad\alpha=1,\ldots,m,\quad I\in\mathbb{M}_{n,k}.
\]
Recall that $J^{\infty}\pi$ is, by definition, an inverse limit of the tower
of projections
\begin{equation}
\xymatrix@C=30pt{M  & E \ar[l]_-{\pi} & \cdots \ar[l]& J^{k}\pi \ar[l]_-{\pi_{k,k-1}} & J^{k+1}\pi \ar[l]_-{\pi_{k+1,k}} & \cdots \ar[l]
}. \label{Jinfty}%
\end{equation}

Now, let $k<\infty$, $\tau_{0}:T_{0}\longrightarrow J^{k}\pi$ be a vector
bundle, $\dim T_{0}=\dim J^{k}\pi+p$, and $(\ldots,x^{i},\ldots,u_{I}^{\alpha
},\ldots,v^{a},\ldots)$ adapted to $\tau_{0}$, local coordinates on $T_{0}$. A
(possibly non-linear) \emph{differential operator of order }$\leq k$ `acting
on local sections of $\pi$, with values in $\tau_{0}$' (in short `from $\pi$
to $\tau_{0}$') is a section $\Phi:J^{k}\pi\longrightarrow T_{0}$ of $\tau
_{0}$. For any local section $\boldsymbol{p}:U\longrightarrow E$ of $\pi$,
$\Phi$ determines an `image' section $\Delta_{\Phi}\boldsymbol{p}:=\Phi\circ
j_{k}\boldsymbol{p}:U\longrightarrow T_{0}$ of the bundle $\underline{\tau
}_{0}:=\pi_{k}\circ\tau_{0}:T_{0}\longrightarrow M$. If $\Phi$ is locally
given by
\begin{equation}
v^{a}=\Phi^{a}(\ldots,x^{i},\ldots,u_{I}^{\beta},\ldots),\quad a=1,\ldots,p,
\label{LocRepF}%
\end{equation}
and $\boldsymbol{p}$ is locally given by (\ref{p}), then $\Delta_{\Phi
}\boldsymbol{p}$ is locally given by
\[
\left\{
\begin{array}
[c]{l}%
u_{I}^{\alpha}=(\tfrac{\partial^{|I|}}{\partial x^{I}}p^{\alpha})(\ldots
,x^{i},\ldots)\\
v^{a}=\Phi^{a}(\ldots,x^{i},\ldots,(\tfrac{\partial^{|J|}}{\partial x^{J}%
}p^{\beta})(\ldots,x^{j},\ldots),\ldots)
\end{array}
\right.  ,
\]
$\alpha=1,\ldots,m,\;I\in\mathbb{M}_{n,k},\;a=1,\ldots,p$. This motivates the
name `differential operator' for $\Phi$. Denote by $\mathrm{diff}_{k}%
(\pi,\tau_{0})$ the set of all differential operators of order $\leq k$ from
$\pi$ to $\tau_{0}$.

For $\Phi\in\mathrm{diff}_{k}(\pi,\tau_{0})$ and $l\leq\infty$ we define the
$l$th prolongation of $\Phi$ as follows. Consider the space $J^{l}%
\underline{\tau}{}_{0}$ of $l$-jets of local sections of $\underline{\tau}%
{}_{0}$, and local jet coordinates $(\ldots,x^{i},\ldots,u_{I,J}^{\alpha
},\ldots,v_{J}^{a},\ldots)$ on $J^{l}\underline{\tau}{}_{0}$, $J\in
\mathbb{M}_{n,l}$. In $J^{l}\underline{\tau}{}_{0}$ consider the submanifold
$T_{0}^{(l)}$ made of jets of local sections of the form $\Delta_{\Psi
}\boldsymbol{p}$, where $\Psi\in\mathrm{diff}_{k}(\pi,\tau_{0})$ and
$\boldsymbol{p}$ is a local section of $\pi$. $T_{0}^{(l)}$ is locally defined
by
\[
u_{I,J}^{\alpha}=u_{I^\prime,J^\prime}^{\alpha},\;\text{whenever}I,I^\prime\in\mathbb{M}%
_{n,k},\;J,J^\prime\in\mathbb{M}_{n,l} \text{ are such that } IJ = I^\prime J^\prime,
\]
$\alpha=1,\ldots,m$. Thus $(\ldots,x^{i},\ldots,u_{I}^{\alpha},\ldots,v_{J}^{a},\ldots)$,
$I\in\mathbb{M}_{n,k+l},\;J\in\mathbb{M}_{n,l}$, are local coordinates on
$T_{0}^{(l)}$. $T_{0}^{(l)}$ projects canonically onto $J^{k+l}\pi$ and the
projection $\tau_{0}^{(l)}:T_{0}^{(l)}\longrightarrow J^{k+l}\pi$ is a vector
bundle. Moreover, coordinates $(\ldots,x^{i},\ldots,u_{I}^{\alpha},\ldots
,v_{J}^{a},\ldots)$ on $T_{0}^{(l)}$ are adapted to $\tau_{0}^{(l)}$. Finally,
define the $l$th prolongation $\Phi^{(l)}:J^{k+l}\pi\longrightarrow
T_{0}^{(l)}$ of $\Phi$ by putting $\Phi^{(l)}([\boldsymbol{p}]_{x}^{k+l}):=[\Delta
_{\Phi}\boldsymbol{p}]_{x}^{l}\in T_{0}^{(l)}$, for all local sections
$\boldsymbol{p}$ of $\pi$ and $x\in M$. Then $\Phi^{(l)}\in\mathrm{diff}%
_{k+l}(\pi,\tau_{0}^{(l)})$.

For $\Phi\in\mathrm{diff}_{k}(\pi,\tau_{0})$ put $\mathscr{E}_{\Phi}%
:=\{\theta\in J^{k}\pi\;|\;\Phi(\theta)=0\}$. $\mathscr{E}_{\Phi}$ is called
the \emph{(system of) PDE(s)} determined by $\Phi$. For $l\leq\infty$ put also
$\mathscr{E}_{\Phi}^{(l)}:=\mathscr{E}_{\Phi^{(l)}}$. $\mathscr{E}_{\Phi
}^{(l)}$ is locally determined by equations
\begin{equation}
(D_{J}\Phi^{a})(\ldots,x^{i},\ldots,u_{I}^{\alpha},\ldots)=0,\quad
a=1,\ldots,p,\;J\in\mathbb{M}_{n,l}, \label{Prolong}%
\end{equation}
where $D_{j_{1}\cdots j_{l}}:=D_{j_{1}}\circ\cdots\circ D_{j_{l}}$ and
$D_{j}:=\partial/\partial x^{j}+u_{Ij}^{\alpha}\partial/\partial u_{I}%
^{\alpha}$ is the $j$\emph{th total derivative, }$j,j_{1},\ldots,j_{l}=1,\ldots,m$.
In the following we put $\partial_{\alpha}^{I}:=\partial/\partial
u_{I}^{\alpha}$ and $\partial_{\alpha}:=\partial/\partial u^{\alpha}$,
$\alpha=1,\ldots,m$, $I\in\mathbb{M}_{n}$.

A local section $\boldsymbol{p}$ of $\pi$ is a \emph{(local) solution of
}$\mathscr{E}_{\Phi}$ iff, by definition, $\Gamma_{\boldsymbol{p}}^{k}%
\subset\mathscr{E}_{\Phi}$ or, which is the same, $\Gamma_{\boldsymbol{p}%
}^{k+l}\subset\mathscr{E}_{\Phi}^{(l)}$ for some $l\leq\infty$. Notice that
$\mathscr{E}_{\Phi}^{(\infty)}\subset J^{\infty}\pi$ is an inverse limit of
the tower of maps
\begin{equation}
\xymatrix@C=45pt{M  &\mathscr{E}_\Phi \ar[l]_-{\pi_k} & \cdots \ar[l]& \mathscr{E}^{(l)}_\Phi \ar[l]_-{\pi_{k+l,k+l-1}} & \mathscr{E}^{(l+1)}_\Phi \ar[l]_-{\pi_{k+l+1,k+l}} & \cdots \ar[l]
} \label{Einfty}%
\end{equation}
and consists of \textquotedblleft formal solutions\textquotedblright\ of
$\mathscr{E}_{\Phi}$, i.e., possibly non-converging Taylor series fulfilling
(\ref{Prolong}) for every $l$. The PDE $\mathscr{E}_{\Phi}$ is called
\emph{formally integrable} iff $\mathscr{E}_{\Phi}^{(l)}\subset J^{k+l}\pi$ is
a (closed) submanifold for any $l<\infty$ and (\ref{Einfty}) is a sequence of
fiber bundles. Let us stress that, basically, all relevant PDEs in
Mathematical Physics are formally integrable and, therefore, in the following,
we will only consider differential operators determining formally integrable PDEs.

$J^{\infty}\pi$ and $\mathscr{E}_{\Phi}^{(\infty)}$ are not finite dimensional
smooth manifolds, in general. However, they are \emph{pro-finite dimensional
smooth manifolds}. We do not give here a complete definition of a pro-finite
dimensional smooth manifold, which would take too much space. Rather, we will
just outline it. Basically, a pro-finite dimensional smooth manifold is a(n
equivalence class of) set(s) $\mathscr{O}$ together with a sequence of smooth
fiber bundles%
\begin{equation}
\xymatrix@C=30pt{\mathscr{O}_0  &\mathscr{O}_1 \ar[l]_-{\mu_{1,0}} & \cdots \ar[l]& \mathscr{O}_{k} \ar[l]_-{\mu_{k,k-1}} & \mathscr{O}_{k+1} \ar[l]_-{\mu_{k+1,k}} & \cdots \ar[l]
} \label{Oinfty}%
\end{equation}
and maps $\mu_{\infty,k}:\mathscr{O}\longrightarrow\mathscr{O}_{k}$, $0\leq
k<\infty$, such that $\mathscr{O}$ (together with the $\mu_{\infty,k}$'s) is
an inverse limit of (\ref{Oinfty}). It is associated to the sequence
(\ref{Oinfty}) a filtration of algebras
\begin{equation}
\xymatrix@C=30pt{C^\infty(\mathscr{O}_{0}) \ar[r]^-{\mu^\ast_{1,0}} & \cdots \ar[r] & C^\infty(\mathscr{O}_{k-1}) \ar[r]^-{\mu_{k,k-1}^\ast} & C^\infty(\mathscr{O}_{k}) \ar[r]^-{\mu_{k+1,k}^\ast} & \cdots 
}. \label{COinfty}%
\end{equation}
We understand the monomorphisms $\mu_{l+1,l}^{\ast}$'s and interpret
(\ref{COinfty}) as a sequence of subalgebras. Similarly, we understand the
$\mu_{\infty,l}$ 's and interpret elements in $C^{\infty}(\mathscr{O}_{k})$ as
functions on $\mathscr{O}$. Put $C^{\infty}(\mathscr{O}):=\bigcup
\nolimits_{l\in\mathbb{N}_{0}}C^{\infty}(\mathscr{O}_{k})$. $C^{\infty
}(\mathscr{O})$ is interpreted as algebra of \emph{smooth functions on
}$\mathscr{O}$. Differential calculus over $\mathscr{O}$ may then be
introduced as \emph{filtered differential calculus over }$C^{\infty
}(\mathscr{O})$ \cite{v01}. Since the main constructions (smooth maps, vector
fields, differential forms, linear jets and differential operators, etc.) of
such calculus do not look very different from the analogous ones in
finite-dimensional differential geometry we will not insist on this and refer
to \cite{v01} for the rigorous definitions and the main results (see
\cite{s89} and \cite{t89,t91} for a sketch of alternative approaches).

Here we just recall the definition of finite dimensional vector bundle over
$\mathscr{O}$. This is, basically, a vector bundle over $\mathscr{O}_{k}$ for
some $k<\infty$, pull-backed to $\mathscr{O}$ via $\mu_{\infty,k}$. In more
details, let $\tau_{0}:T_{0}\longrightarrow\mathscr{O}_{k}$ be a (finite
dimensional) vector bundle, $k<\infty$. For $l\geq0$ let $\tau_{l}%
:=\mu_{k+l,k}^{\circ}(\tau_{0}):T_{l}:=\mu_{k+l,k}^{\circ}(T_{0}%
)\longrightarrow\mathscr{O}_{k+l}$ be the induced (by $\tau_{0}$ via
$\mu_{k+l,k}$) vector bundle and $\nu_{l+1,l}:T_{l+1}\longrightarrow T_{l}$
the canonical projection. Denote by $T$ the pro-finite dimensional smooth
manifold determined by the sequence of fiber bundles
\begin{equation}
\xymatrix@C=40pt{{T}_0  &{T}_1 \ar[l]_-{\nu_{1,0}} & \cdots \ar[l]& {T}_{l} \ar[l]_-{\nu_{l,l-1}} & {T}_{l+1} \ar[l]_-{\nu_{l+1,l}} & \cdots \ar[l]
}. \label{Tinfty}%
\end{equation}
The maps $\tau_{l}:T_{l}\longrightarrow\mathscr{O}_{l+k}$, $l\geq0$, determine
a smooth map $\tau:T\longrightarrow\mathscr{O}$. Any such map is, by
definition, \emph{a (finite-dimensional) vector bundle over} $\mathscr{O}$.
Notice that it is associated to the sequence (\ref{Tinfty}) of vector bundle
morphisms a filtration of vector spaces
\[
\xymatrix@C=50pt{\Gamma(\tau_{0}) \ar[r]^-{\mu^\circ_{k+1,k}} & \cdots \ar[r] & \Gamma(\tau_{l-1}) \ar[r]^-{\mu_{k+l,k+l-1}^\circ} & \Gamma(\tau_{l}) \ar[r]^-{\mu_{k+l+1,k+l}^\circ} & \cdots 
}.
\]
We understand the monomorphisms $\mu_{k+l+1,k+l}^{\circ}$'s and interpret
(\ref{COinfty}) as a sequence of vector subspaces. Similarly, we understand
the $\mu_{\infty,k+l}$'s and interpret elements in $\Gamma(\tau_{l})$ as
functions $\mathscr{O}\longrightarrow T$. Put $\Gamma(\tau):=\bigcup
\nolimits_{l\in\mathbb{N}_{0}}\Gamma(\tau_{l})$. $\Gamma(\tau)$ is naturally a
$C^{\infty}(\mathscr{O})$-module and it is interpreted as the module of
\emph{smooth sections of }$\tau$.

As an example, let $\mathscr{O}=J^{\infty}\pi$, $\tau_{0}:T_{0}\longrightarrow
J^{k}\pi$ be a vector bundle for some $k<\infty$ and $\tau:=\pi_{\infty
,k}^{\circ}(\tau_{0}):T:=\pi_{\infty,k}^{\circ}(T_{0})\longrightarrow
J^{\infty}\pi$. Since $\Gamma(\tau_{l})=\mathrm{diff}_{k+l}(\pi,\tau_{l})$ for
any $l$, we have the filtration $\mathrm{diff}_{k}(\pi,\tau_{0})\subset
\mathrm{diff}_{k+1}(\pi,\tau_{1})\subset\cdots\subset\mathrm{diff}_{k+l}%
(\pi,\tau_{l})\subset\cdots$. Put $\mathrm{diff}(\pi,\tau):=\bigcup
_{l\in\mathbb{N}_{0}}\mathrm{diff}_{k+l}(\pi,\tau_{l})=\Gamma(\tau)$. Elements
in $\mathrm{diff}(\pi,\tau)$ are called \emph{differential operators} `acting
on local sections of $\pi$, with values in $\tau_{0}$' (in short `from $\pi$
to $\tau_{0}$'). They are nothing but sections of the vector bundle
$\tau:T\longrightarrow J^{\infty}\pi$.

An important technical advantage of formally integrable PDEs is the following.
Let $\mathscr{E}\subset J^{\infty}(\pi)$ be the $\infty$th prolongation of a
formally integrable PDE, $\tau:T\longrightarrow J^{\infty}(\pi)$ a vector
bundle and $\tau|_{\mathscr{E}}:T|_{\mathscr{E}}\longrightarrow\mathscr{E}$
its restriction to $\mathscr{E}$. Then for any section $s\in\Gamma
(\tau|_{\mathscr{E}})$ there exists a section $\widetilde{s}\in\Gamma(\tau)$
such that $s=\tilde{s}|_{\mathscr{E}}$. In the following we will often use
this property without further comments.

Finally, let us mention here that a vector field on an pro-finite dimensional
manifold does not generate a flow in general (see, for instance, \cite{c91}).

\subsection{The Cartan Distribution and the $\mathscr{C}$-Spectral
Sequence\label{Sec2}}

Let $\pi:E\longrightarrow M$ and $\tau:T\longrightarrow J^{\infty}\pi$ be as
in the previous section and $\Phi\in\mathrm{diff}(\pi,\tau)$. In the following
we will simply write $J^{\infty}$ for $J^{\infty}\pi$ and $\mathscr{E}$ for
$\mathscr{E}_{\Phi}^{(\infty)}$. $i_{\mathscr{E}}:\mathscr{E}\hookrightarrow
J^{\infty}$ will denote the inclusion. Notice that for $\Phi=0$,
$\mathscr{E}=\mathscr{E}_{\Phi}^{(\infty)}=J^{\infty}$.

Recall that $J^{\infty}$ is endowed with the Cartan distribution $\mathscr{C}$
which is defined as follows:
\[
\mathscr{C}:J^{\infty}\ni\theta\longmapsto\mathscr{C}_{\theta}\subset
T_{\theta}J^{\infty},
\]
where $\mathscr{C}_{\theta}:=T_{\theta}\Gamma_{\boldsymbol{p}}^{\infty}$ for
$\theta=[\boldsymbol{p}]_{x}^{\infty}$, $x\in M$. Denote by
$\mathscr{C}\mathrm{D}(J^{\infty})\subset\mathrm{D}(J^{\infty})$ the
$C^{\infty}(J^{\infty})$-submodule made of vector fields in the Cartan
distribution, i.e., vector fields $X\in\mathrm{D}(J^{\infty})$ such that
$X_{\theta}\in\mathscr{C}_{\theta}$ for all $\theta\in J^{\infty}$. The Cartan
distribution is $n$-dimensional, it is locally spanned by total derivatives
$\ldots,D_{i},\ldots$ and it is involutive, i.e., $[X,Y]\in
\mathscr{C}\mathrm{D}(J^{\infty})$ for all $X,Y\in\mathscr{C}\mathrm{D}%
(J^{\infty})$. Moreover, $n$-dimensional integral submanifolds $L\subset
J^{\infty}$ of $\mathscr{C}$ are of the form $L=\Gamma_{\boldsymbol{p}%
}^{\infty}$ for some local section $\boldsymbol{p}$ of $\pi$.

Let $\mathscr{E}\subset J^{\infty}$ be as above. The
Cartan distribution $\mathscr{C}$ restricts to $\mathscr{E}$ in the sense that
$\mathscr{C}_{\theta}\subset T_{\theta}\mathscr{E}$ for any $\theta
\in\mathscr{E}$. Abusing the notation we still denote by $\mathscr{C}$ the
restricted to $\mathscr{E}$ distribution and call it the \emph{Cartan
distribution of }$\mathscr{E}$. Also we denote by $\mathscr{C}\mathrm{D}%
(\mathscr{E})\subset\mathrm{D}(\mathscr{E})$ the $C^{\infty}(\mathscr{E})$%
-submodule made of vector fields in $\mathscr{C}$. Elements in
$\mathscr{C}\mathrm{D}(\mathscr{E})$ are called \emph{horizontal vector
fields}. In particular, total derivatives restrict to $\mathscr{E}$, i.e.,
there are unique local vector fields $\ldots,D_{i}^{\mathscr{E}},\ldots$ on
$\mathscr{E}$ such that $i_{\mathscr{E}}^{\ast}\circ D_{i}=D_{i}%
^{\mathscr{E}}\circ i_{\mathscr{E}}^{\ast}$, $i=1,\ldots,n$. Again
$\mathscr{C}$ is locally spanned by vector fields $\ldots,D_{i}^{\mathscr{E}}%
,\ldots$, it is involutive and $n$-dimensional integral submanifolds of it
are graphs $\Gamma_{\boldsymbol{p}}^{\infty}$ of infinite jet prolongations of
local solutions $\boldsymbol{p}$ of $\mathscr{E}_{\Phi}$.

A spectral sequence is naturally associated to an involutive distribution
and, in particular, to the Cartan distribution on (the infinite prolongation
of) a PDE as follows. Denote by $\mathscr{C}\Lambda(\mathscr{E})\subset
\Lambda(\mathscr{E})$ the subset made of differential forms $\omega$ such
that
\[
\omega(X_{1},\ldots,X_{k})=0\quad\text{for all }X_{1},\ldots,X_{k}%
\in\mathscr{C}\mathrm{D}(\mathscr{E})\text{,}%
\]
where $k$ is the degree of $\omega$. $\mathscr{C}\Lambda(\mathscr{E})$ is a
differential ideal in $\Lambda(\mathscr{E})$. Namely, it is an algebraic ideal
and, moreover, it is differentially closed, i.e., $d\omega\in
\mathscr{C}\Lambda(\mathscr{E})$ for any $\omega\in\mathscr{C}\Lambda
(\mathscr{E})$. $\mathscr{C}\Lambda(\mathscr{E})$ is called the \emph{Cartan
ideal} of $\mathscr{E}$. For any $p\in\mathbb{N}$, denote by $\mathscr{C}^{p}%
\Lambda(\mathscr{E})$ the $p$th exterior power of $\mathscr{C}\Lambda
(\mathscr{E})$. Thus, the sequence
\[
\Lambda(\mathscr{E})\supset\mathscr{C}\Lambda(\mathscr{E})\supset
\mathscr{C}^{2}\Lambda(\mathscr{E})\supset\cdots\supset\mathscr{C}^{p}%
\Lambda(\mathscr{E})\supset\cdots
\]
is a filtration of the de Rham complex $(\Lambda(\mathscr{E}),d)$ of
$\mathscr{E}$. The associated spectral sequence is denoted by
$\mathscr{C}E(\mathscr{E})=\{(\mathscr{C}E_{r}^{p,q}(\mathscr{E}),d_{r}%
^{p,q})\}_{r}^{p,q}$ and called the $\mathscr{C}$-spectral sequence of
$\mathscr{E}$ \cite{v84}. It is regular and converges to de Rham cohomologies
of $\mathscr{E}$.

The first column of the $0$th term of $\mathscr{C}E(\mathscr{E})$,
\[
\xymatrix{0 \ar[r] & \mathscr{C}E_0^{0,0}(\mathscr{E}) \ar[r]^-{d_0^{0,0}} & \mathscr{C}E_0^{0,1}(\mathscr{E}) \ar[r]^-{d_0^{0,1}} & \cdots \ar[r] & \mathscr{C}E_0^{0,q}(\mathscr{E}) \ar[r]^-{d_0^{0,q}} & \cdots},
\]
is, by definition, the quotient complex $\Lambda
(\mathscr{E})/\mathscr{C}\Lambda(\mathscr{E})$, which is also denoted by
\[
\xymatrix{0 \ar[r] &  C^\infty(\mathscr{E}) \ar[r]^-{\overline{d}} & \overline{\Lambda}{}^1(\mathscr{E}) \ar[r]^-{\overline{d}} & \cdots \ar[r] & \overline{\Lambda}{}^q(\mathscr{E}) \ar[r]^-{\overline{d}} & \cdots},
\]
and called \emph{the horizontal de Rham complex of }$\mathscr{E}$. Its
cohomology algebra $\mathscr{C}E_{1}^{0,\bullet}(\mathscr{E})$ is denoted by $\overline
{H}(\mathscr{E})$, and called \emph{horizontal de Rham cohomology algebra of}
$\mathscr{E}$. Recall, in particular, that $\overline{d}$-closed elements in
$\overline{\Lambda}{}^{n-1}(\mathscr{E})$ are called \emph{conserved currents}
and cohomology classes in $\overline{H}{}^{n-1}(\mathscr{E})$
\emph{conservation laws} of the PDE $\mathscr{E}_{\Phi}$.

In the following we will denote by $\mathscr{C}\Lambda^{k}(\mathscr{E})$
(resp.~$\mathscr{C}^{p}\Lambda^{k}(\mathscr{E})$, $\overline{\Lambda}{}%
^{k}(\mathscr{E})$, $\overline{H}{}^{k}(\mathscr{E})$) the $k$th homogeneous
component of $\mathscr{C}\Lambda(\mathscr{E})$ (resp. $\mathscr{C}^{p}%
\Lambda(\mathscr{E})$, $\overline{\Lambda}{}(\mathscr{E})$, $\overline
{H}(\mathscr{E})$), $k\geq0$, and by $\mathscr{C}^{\bullet}\Lambda
(\mathscr{E}):=\bigoplus_{p}\mathscr{C}^{p}\Lambda^{p}(\mathscr{E})\subset
\Lambda(\mathscr{E})$ the $C^{\infty}(\mathscr{E})$-subalgebra generated by
$\mathscr{C}\Lambda^{1}(\mathscr{E})$. Notice that $\mathscr{C}^{p}%
\Lambda(\mathscr{E})$ is generated by $\mathscr{C}^{p}\Lambda^{p}%
(\mathscr{E})$ as an ideal, $p>0$.

The $\mathscr{C}$-spectral sequence $\mathscr{C}E(\mathscr{E})$ contains very
relevant \textquotedblleft invariants\textquotedblright\ of the PDE
$\mathscr{E}_{\Phi}$ (see, for instance, \cite{b...99,kv98}). Moreover, it
formalizes in a coordinate-free manner variational calculus (on local sections
of $\pi$) constrained by $\mathscr{E}_{\Phi}$ \cite{v84}. Therefore, it is a
most fundamental construction in the geometric theory of differential
equations. Finally, it is a very general construction. For instance, it may be
defined exactly in the same way when $\mathscr{E}$ is the infinite
prolongation of a system of PDEs \textquotedblleft imposed on general
$n$-dimensional submanifolds of $E$\textquotedblright. However, in the
present case, the fibered structure $\pi_{\infty}|_{\mathscr{E}}%
:\mathscr{E}\longrightarrow M$ of $\mathscr{E}$ allows a more simple
description (which is, in the general case, valid only locally), \emph{the
variational bi-complex }\cite{v84}, which we briefly recall in the following.

The Cartan distribution and the fibered structure $\pi_{\infty}|_{\mathscr{E}}%
:\mathscr{E}\longrightarrow M$ of $\mathscr{E}$ determine a splitting of the
tangent bundle $T\mathscr{E}\longrightarrow\mathscr{E}$ into the Cartan or
horizontal part $\mathscr{C}$ and the vertical (with respect to $\pi_{\infty}%
$) part $V\pi_{\infty}|_{\mathscr{E}}$. Accordingly, $\mathrm{D}(\mathscr{E})$
splits into a direct sum: $\mathrm{D}(\mathscr{E})=\mathscr{C}\mathrm{D}%
(\mathscr{E})\oplus V\mathrm{D}(\mathscr{E})$, $V\mathrm{D}%
(\mathscr{E})\subset\mathrm{D}(\mathscr{E})$ being the $C^{\infty
}(\mathscr{E})$-submodule made of $\pi_{\infty}$\emph{-vertical vector
fields}, i.e., vector fields $Y\in\mathrm{D}(\mathscr{E})$ such that
$Y\circ\pi_{\infty}^{\ast}=0$. In particular, $V\mathrm{D}(J^{\infty})$ is
locally (formally) generated by vector fields $\ldots,\partial_{\alpha}^{I},\ldots$. Dually,
$\Lambda^{1}(\mathscr{E})$ splits into the direct sum
\begin{equation}
\Lambda^{1}(\mathscr{E})=\mathscr{C}\Lambda^{1}(\mathscr{E})\oplus
\overline{\Lambda}{}^{1}(\mathscr{E}); \label{SplitLambda1}%
\end{equation}
here and in what follows $\overline{\Lambda}{}^{1}(\mathscr{E})$ is identified
with the $C^{\infty}(\mathscr{E})$-submodule in $\Lambda^{1}(\mathscr{E})$
generated by $\pi_{\infty}^{\ast}(\Lambda^{1}(M))$. In particular,
$\mathscr{C}\Lambda^{1}(J^{\infty})$ is locally generated by forms
$\ldots,\omega_{I}^{\alpha}:=du_{I}^{\alpha}-u_{Ii}^{\alpha}dx^{i},\ldots$ and
$\overline{\Lambda}{}^{1}(J^{\infty})$ is locally generated by forms
$\ldots,dx^{i},\ldots$. Similarly, $\mathscr{C}\Lambda^{1}(\mathscr{E})$ is
locally generated by forms $\ldots,i_{\mathscr{E}}^{\ast}(\omega_{I}^{\alpha
}),\ldots$ and $\overline{\Lambda}{}^{1}(\mathscr{E})$ is locally generated by
forms $\ldots,i_{\mathscr{E}}^{\ast}(dx^{i}),\ldots$.

In view of splitting (\ref{SplitLambda1}) $\Lambda(\mathscr{E})$ factorizes as
$\Lambda(\mathscr{E})\simeq\mathscr{C}^{\bullet}\Lambda(\mathscr{E})\otimes
\overline{\Lambda}{}(\mathscr{E})$ (here and in what follows tensor products
will be always over $C^{\infty}(\mathscr{E})$, or $C^{\infty}(J^{\infty})$ for
$\Phi=0$). In particular, there are projections $\mathfrak{p}_{p,q}%
:\Lambda(\mathscr{E})\longrightarrow\mathscr{C}^{p}\Lambda^{p}%
(\mathscr{E})\otimes\overline{\Lambda}{}^{q}(\mathscr{E})$ for any
$p,q\in\mathbb{N}_{0}$. Correspondingly, the de Rham complex of $\mathscr{E}$,
$(\Lambda(\mathscr{E}),d)$, splits in a bi-complex $(\mathscr{C}^{\bullet
}\Lambda(\mathscr{E})\otimes\overline{\Lambda}{}(\mathscr{E}),\overline
{d},d^{V})$ (in the following diagram we drop for simplicity the postfix
$(\mathscr{E})$),
\begin{equation}%
\begin{array}
[c]{c}%
\xymatrix{& \cdots & \cdots & \cdots & \cdots & \cdots \\
0 \ar[r]& \mathscr{C}^{p+1}\Lambda^{p+1} \ar[r]^-{\overline{d}} \ar[u]^-{d^V}& \cdots \ar[r] & \mathscr{C}^{p+1}\Lambda^{p+1}\otimes \overline{\Lambda}{}^q \ar[u]^-{d^V} \ar[r]^-{\overline{d}} &\mathscr{C}^{p+1}\Lambda^{p+1} \otimes \overline{\Lambda}{}^{q+1}  \ar[u]^-{d^V} \ar[r]^-{\overline{d}} & \cdots \\
0 \ar[r] & \mathscr{C}^{p}\Lambda^{p} \ar[r]^-{\overline{d}} \ar[u]^-{d^V} & \cdots \ar[r] & \mathscr{C}^{p}\Lambda^{p}\otimes \overline{\Lambda}{}^{q}  \ar[u]^-{d^V} \ar[r]^-{\overline{d}} &\mathscr{C}^{p}\Lambda^{p}\otimes \overline{\Lambda}{}^{q+1} \ar[u]^-{d^V} \ar[r]^-{\overline{d}} & \cdots \\
 & \cdots \ar[u] & \cdots & \cdots \ar[u] & \cdots \ar[u] & \cdots \\
0 \ar[r] & C^\infty \ar[r]^-{\overline{d}} \ar[u]^-{d^V} & \cdots \ar[r] & \overline{\Lambda}{}^q \ar[r]^-{\overline{d}} \ar[u]^-{d^V} & 
\overline{\Lambda}{}^{q+1} \ar[r]^-{\overline{d}} \ar[u]^-{d^V} & \cdots \\
& 0 \ar[u] & & 0 \ar[u] & 0 \ar[u] &}
\end{array}
, \label{VarBicomp}%
\end{equation}
defined by
\[
\overline{d}(\omega\otimes\overline{\sigma}):=(\mathfrak{p}_{p,q+1}\circ
d)(\omega\wedge\overline{\sigma})\text{\quad and}\quad d^{V}(\omega
\otimes\overline{\sigma}):=(\mathfrak{p}_{p+1,q}\circ d)(\omega\wedge
\overline{\sigma}),
\]
where $\omega\in\mathscr{C}^{p}\Lambda^{p}(\mathscr{E})$ and $\overline
{\sigma}\in\overline{\Lambda}{}^{q}(\mathscr{E})$, $p,q\in$ $\mathbb{N}_{0}$.
$\overline{d}$ and $d^{V}$ are called the \emph{horizontal} and the
\emph{vertical de Rham differential}, respectively, and (\ref{VarBicomp}) is
called the \emph{variational bi-complex}. In the following we will often
understand isomorphism $\Lambda(\mathscr{E})\simeq\mathscr{C}^{\bullet}%
\Lambda(\mathscr{E})\otimes\overline{\Lambda}{}(\mathscr{E})$.

As a bi-complex (\ref{VarBicomp}) determines two spectral sequences. One of
them is the $\mathscr{C}$-spectral sequence while the other is the
Leray-Serre spectral sequence of the fibration $\pi_{\infty}|_{\mathscr{E}}%
:\mathscr{E}\longrightarrow M$ \cite{mvv??}. In particular, for any $p$, there
is a canonical isomorphisms of complexes
\begin{equation}
(\mathscr{C}E_{0}^{p,\bullet}(\mathscr{E}),d_{0}^{p,\bullet})\simeq
(\mathscr{C}^{p}\Lambda^{p}(\mathscr{E})\otimes\overline{\Lambda}%
{}(\mathscr{E}),\overline{d}), \label{HorDRComp}%
\end{equation}
and the differential $d_{1}^{p,\bullet}:\mathscr{C}E_{1}^{p,\bullet
}(\mathscr{E})\longrightarrow\mathscr{C}E_{1}^{p+1,\bullet}(\mathscr{E})$ is isomorphic to the map
induced by $d^{V}$ in the cohomology $H(\mathscr{C}^{p}\Lambda^{p}%
(\mathscr{E})\otimes\overline{\Lambda}{}(\mathscr{E}),\overline{d})$.

Notice that the embedding $i_{\mathscr{E}}:\mathscr{E}\hookrightarrow
J^{\infty}$ of the infinite prolongation $\mathscr{E}$ of a PDE determines via
pull-back both a morphism of spectral sequences and a morphism of
bi-complexes that, abusing the notation, we denote by the same symbol
\begin{gather*}
i_{\mathscr{E}}^{\ast}:\{(\mathscr{C}E_{r}^{\bullet,\bullet}(J^{\infty}%
),d_{r}^{\bullet,\bullet})\}\longrightarrow\{(\mathscr{C}E_{r}^{\bullet
,\bullet}(\mathscr{E}),d_{r}^{\bullet,\bullet})\},\\
i_{\mathscr{E}}^{\ast}:(\mathscr{C}^{\bullet}\Lambda(J^{\infty})\otimes
\overline{\Lambda}{}(J^{\infty}),\overline{d},d^{V})\longrightarrow
(\mathscr{C}^{\bullet}\Lambda(\mathscr{E})\otimes\overline{\Lambda}%
{}(\mathscr{E}),\overline{d},d^{V}).
\end{gather*}
\subsection{Higher Symmetries of PDEs\label{Sec3}}

Denote by $\mathrm{D}_{\mathscr{C}}(\mathscr{E})\subset\mathrm{D}%
(\mathscr{E})$ the subset made of vector fields preserving the Cartan
distribution, i.e., vector fields $X$ such that $[X,Y]\in\mathscr{C}\mathrm{D}%
(\mathscr{E})$ for any $Y\in\mathscr{C}\mathrm{D}(\mathscr{E})$.
$\mathrm{D}_{\mathscr{C}}(\mathscr{E})$ is clearly a Lie subalgebra in
$\mathrm{D}(\mathscr{E})$. Elements in $\mathrm{D}_{\mathscr{C}}(\mathscr{E})$
are called \emph{(infinitesimal) symmetries of }$\mathscr{E}_{\Phi}$. The
theory of infinitesimal symmetries of PDEs is fundamental in many respects
\cite{b...99}. Notice that, since the Cartan distribution is involutive, then
$\mathscr{C}\mathrm{D}(\mathscr{E})\subset\mathrm{D}_{\mathscr{C}}%
(\mathscr{E})$ and it is an ideal in $\mathrm{D}_{\mathscr{C}}(\mathscr{E})$.
Elements in $\mathscr{C}\mathrm{D}(\mathscr{E})$ are called \emph{trivial
symmetries of }$\mathscr{E}_{\Phi}$, in that horizontal vector fields
\textquotedblleft are symmetries of every PDE\textquotedblright. The quotient
Lie algebra $\mathrm{Sym}(\mathscr{E}):=\mathrm{D}_{\mathscr{C}}%
(\mathscr{E})/\mathscr{C}\mathrm{D}(\mathscr{E})$ is called the algebra of
\emph{non-trivial higher symmetries of }$\mathscr{E}_{\Phi}$. Clearly, every
equivalence class $\boldsymbol{X}=X+\mathscr{C}\mathrm{D}(\mathscr{E})\in
\mathrm{Sym}(\mathscr{E})$, $X\in\mathrm{D}_{\mathscr{C}}(\mathscr{E})$, has
got one and only one vertical representative $X^{V}\in V\mathrm{D}%
(\mathscr{E})$. Any vertical element in $\mathrm{D}_{\mathscr{C}}%
(\mathscr{E})$ is called an \emph{evolutionary vector field}. Thus
$\mathrm{Sym}(\mathscr{E})$ is isomorphic to the Lie algebra $V\mathrm{D}%
_{\mathscr{C}}(\mathscr{E})$ of evolutionary vector fields.

In order to effectively describe $V\mathrm{D}_{\mathscr{C}}(\mathscr{E})$ and,
therefore, $\mathrm{Sym}(\mathscr{E})$ let us first consider the case
$\mathscr{E}=J^{\infty}$. It is easy to prove that any evolutionary vector field
$Y\in V\mathrm{D}_{\mathscr{C}}(J^{\infty})$ is determined by its restriction
to $C^{\infty}(E)\subset C^{\infty}(J^{\infty})$. Moreover, every vertical
vector field $\chi:C^{\infty}(E)\longrightarrow C^{\infty}(J^{\infty})$ along
$\pi_{\infty,0}:J^{\infty}\longrightarrow E$ ($\chi$ is vertical if $\chi
\circ\pi^{\ast}=0$) extends to a unique evolutionary vector field $\rE_{\chi
}\in V\mathrm{D}_{\mathscr{C}}(J^{\infty})$. We conclude that $V\mathrm{D}%
_{\mathscr{C}}(J^{\infty})$ is in one to one correspondence with the
$C^{\infty}(J^{\infty})$-module $\varkappa$ of vector fields along
$\pi_{\infty,0}$ that are vertical with respect to $\pi$ or, which is the same, the module of sections of the induced
vector bundle $\pi_{\infty,0}^{\circ}(\nu\pi):\pi_{\infty,0}^{\circ}%
(V\pi)\longrightarrow J^{\infty}$. Elements in $\varkappa$ are called
\emph{generating sections of higher symmetries} of $\pi$.

Let us now come to the general case when $\mathscr{E}$ is any. First of all
consider the $C^{\infty}(\mathscr{E})$-module $\varkappa|_{\mathscr{E}}$ of
vertical vector fields $\chi:C^{\infty}(E)\longrightarrow C^{\infty
}(\mathscr{E})$ along $\pi_{\infty,0}|_{\mathscr{E}}%
:\mathscr{E}\longrightarrow E$ or, which is the same, the module of sections
of the induced vector bundle $\pi_{\infty,0}|_{\mathscr{E}}^{\circ}(\nu
\pi):\pi_{\infty,0}|_{\mathscr{E}}^{\circ}(V\pi)\longrightarrow\mathscr{E}$.
Elements in $\varkappa|_{\mathscr{E}}$ are called generating sections of
higher symmetries of $\mathscr{E}$. Similarly as to above, a generating
section $\chi\in\varkappa|_{\mathscr{E}}$ extends to a unique vertical vector
field $\rE_{\chi}:C^{\infty}(J^{\infty})\longrightarrow C^{\infty
}(\mathscr{E})$ along the inclusion $i_{\mathscr{E}}%
:\mathscr{E}\hookrightarrow J^{\infty}$. If $\chi$ is locally given by
$\chi=\chi^{\alpha}\partial_{\alpha}$, where $\ldots,\chi^{\alpha},\ldots$ are
local functions on $\mathscr{E}$, then $\rE_{\chi}$ is locally given by
$\rE_{\chi}=D_{I}^{\mathscr{E}}\chi^{\alpha}\partial_{\alpha}^{I}%
|_{\mathscr{E}}$. However, in general $\rE_{\chi}$ is not tangent to
$\mathscr{E}$ and, therefore, is not in $V\mathrm{D}_{\mathscr{C}}%
(\mathscr{E})$. Generating sections $\chi$ such that $\rE_{\chi}\in
V\mathrm{D}_{\mathscr{C}}(\mathscr{E})$ are the ones in the kernel of a
suitable differential operator: the so-called \emph{universal linearization
of }$\mathscr{E}$, which we now define (notice that to the author's knowledge
the following definition never appeared in the literature before in the
general form presented here - see also \cite{t91}).

Let $\tau$, $\Phi$ and $\mathscr{E}$ be as in the previous section, and put
$\underline{\tau}:=\pi_{\infty}\circ\tau:T\longrightarrow M$. Since $\tau$ is
a vector bundle, $V\tau\longrightarrow T$ is naturally isomorphic to the
induced bundle $\tau^{\circ}(\tau):\tau^{\circ}(T)\longrightarrow T$,
$\tau^{\circ}(\tau)$ being the (restriction to $\tau^{\circ}(T)\subset T\times
T$ of the) projection ($T\times T\longrightarrow T$) onto the first factor.
Denote by $\rho_{2}:\tau^{\circ}(T)\longrightarrow T$ the projection onto the
second factor and by $\rho_{2}^{\prime}:V\tau\longrightarrow T$ the map
induced by $\rho_{2}$ via the isomorphism $V\tau\simeq\tau^{\circ}(T)$.
Consider the vertical tangent map $V\Phi:V\pi_{\infty}\longrightarrow
V\underline{\tau}$. Put $o_{\mathscr{E}}:=o\circ i_{\mathscr{E}}%
:\mathscr{E}\longrightarrow T$ and notice, preliminarily, that
$o_{\mathscr{E}}=\Phi\circ i_{\mathscr{E}}$. The short exact sequence of
induced bundles $0\longrightarrow o_{\mathscr{E}}^{\circ}(V\tau
)\longrightarrow o_{\mathscr{E}}^{\circ}(V\underline{\tau})\longrightarrow
V\pi_{\infty}|_{\mathscr{E}}\longrightarrow0$ splits naturally via the map
$Vo|_{\mathscr{E}}:V\pi_{\infty}|_{\mathscr{E}}\longrightarrow o_{\mathscr{E}}%
^{\circ}(V\underline{\tau})$ well defined by putting $Vo|_{\mathscr{E}}%
(\theta,\xi):=(\theta,Vo(\xi))$, $(\theta,\xi)\in V\pi_{\infty}|_{\mathscr{E}}%
$. In particular, there is a canonical projection $V_{\Phi}:o_{\mathscr{E}}%
^{\circ}(V\underline{\tau})\longrightarrow o_{\mathscr{E}}^{\circ}(V\tau).$
Define a map
\[
L_{\Phi}:V\pi_{\infty}|_{\mathscr{E}}\longrightarrow T|_{\mathscr{E}}%
\]
by putting $L_{\Phi}(\xi):=(\theta,\rho_{2}^{\prime}(V))$, where
$(\theta,V):=V_{\Phi}(\theta,V\Phi(\xi))\in o_{\mathscr{E}}^{\circ}(V\tau)$,
for all $\xi\in V_{\theta}\pi_{\infty}$, $\theta\in \mathscr{E}$. $L_{\Phi}$ is
a morphism of vector bundles. For any $\chi\in\varkappa|_{\mathscr{E}}$ let
$\ell_{\Phi}\chi\in\Gamma(\tau|_{\mathscr{E}})$ be defined by putting
$(\ell_{\Phi}\chi)_{\theta}:=L_{\Phi}((\rE_{\chi})_{\theta})$, $\theta
\in\mathscr{E}$. $\ell_{\Phi}:\varkappa|_{\mathscr{E}}\longrightarrow
\Gamma(\tau|_{\mathscr{E}})$ is a well defined linear differential operator
called the \emph{universal linearization of }$\Phi$.

Let us describe $\ell_{\Phi}$ locally. Let $(\ldots,x^{i},\ldots,u_{I}%
^{\alpha},\ldots)$ be local jet coordinates on $J^{\infty}$, $(\ldots
,x^{i},\ldots,u_{I}^{\alpha},\ldots,v^{a},\ldots)$ adapted to $\tau$ local
coordinates on $T$, and $(\ldots,e_{a},\ldots)$ the local basis of
$\Gamma(\tau|_{\mathscr{E}})$ associated to them. If $\Phi$ has local
representation (\ref{LocRepF}), $\ldots,\Phi^{a},\ldots$ being local functions
on $J^{\infty}$, and $\chi=\chi^{\alpha}\partial_{\alpha}$ locally, then
\[
\ell_{\Phi}\chi=e_{a}(\partial_{\alpha}^{I}\Phi^{a})|_{\mathscr{E}}%
D_{I}^{\mathscr{E}}\chi^{\alpha}%
\]
locally.

Now let $\chi\in\varkappa|_{\mathscr{E}}$. It is easy to see that if
$\ell_{\Phi}\chi=0$ then $\rE_{\chi}$ is tangent to $\mathscr{E}$ and,
therefore, it is in $V\mathrm{D}_{\mathscr{C}}(\mathscr{E})$. Vice versa, any
symmetry $Y\in V\mathrm{D}_{\mathscr{C}}(\mathscr{E})$ is of the form
$\rE_{\chi}$ for a unique $\chi\in\varkappa|_{\mathscr{E}}$ such that
$\ell_{\Phi}\chi=0$. We conclude that $\mathrm{Sym}(\mathscr{E})$ is in one to
one correspondence with $\ker\ell_{\Phi}$. In particular, $\ker\ell_{\Phi}$
inherits from $\mathrm{Sym}(\mathscr{E})$ the Lie algebra structure. The
corresponding bracket is denoted by $\{\cdot,\cdot\}$ and called the
\emph{higher Jacobi bracket} of the equation $\mathscr{E}_{\Phi}$.

Finally, notice that, for any $\chi\in\ker\ell_{\Phi}$, the `insertion of' and
the `Lie derivative along' $\rE_{\chi}\in V\mathrm{D}(\mathscr{E})$ commute
with the horizontal de Rham differential $\overline{d}:\Lambda
(\mathscr{E})\longrightarrow\Lambda(\mathscr{E})$, i.e.,
\begin{equation}
i_{\rE_{\chi}}\circ\overline{d}+\overline{d}\circ i_{\rE_{\chi}}%
=\mathcal{L}_{\rE_{\chi}}\circ\overline{d}-\overline{d}\circ\mathcal{L}%
_{\rE_{\chi}}=0. \label{Comm}%
\end{equation}
In their turn Identities (\ref{Comm}) imply
\[
i_{\rE_{\chi}}\circ d^{V}+d^{V}\circ i_{\rE_{\chi}}=\mathcal{L}_{\rE_{\chi}%
},\quad\mathcal{L}_{\rE_{\chi}}\circ d^{V}-d^{V}\circ\mathcal{L}_{\rE_{\chi}%
}=0,
\]
$d^{V}:\Lambda(\mathscr{E})\longrightarrow\Lambda(\mathscr{E})$ being the
vertical de Rham differential.

\subsection{Secondary Differential Forms and Vector Fields\label{Sec4}}

Let $\mathscr{E}$ be as in the previous section. As noticed above,
$n$-dimensional integral submanifolds of the Cartan distribution
$\mathscr{C}$ over $\mathscr{E}$ are in one-to-one correspondence with local
solutions of $\mathscr{E}_{\Phi}$. Thus, informally speaking, the pair
$(\mathscr{E},\mathscr{C})$ encodes all the information about the
\textquotedblleft functional space $\boldsymbol{M}$ of
solutions\textquotedblright\ of $\mathscr{E}_{\Phi}$ (in the following we will
in fact identify $(\mathscr{E},\mathscr{C})$ with $\boldsymbol{M}$). For
instance, \textquotedblleft local functional calculus\textquotedblright\ over
such functional space may be formalized geometrically (and homologically) by
using $(\mathscr{E},\mathscr{C})$ as a starting point and the associated
$\mathscr{C}$-spectral sequence as the main structure. Such formalization has
been named \emph{secondary calculus} \cite{v01} by its discoverer, A. M.
Vinogradov, and its simplest constructions will be briefly reviewed in this section.

Suppose temporarily that $M$ is a compact, orientable and oriented manifold
without boundary. Then an element $\boldsymbol{S}=[\mathscr{L}]\in\overline
{H}{}^{n}(\mathscr{E})=\mathscr{C}E_{1}^{0,n}(\mathscr{E})$,
$\mathscr{L}\in\overline{\Lambda}{}^{n}(\mathscr{E})$, identifies with the
(local) action functional
\[
\boldsymbol{M}\ni\boldsymbol{p}\longmapsto\boldsymbol{S}(\boldsymbol{p}%
):=\int_{M}(j^{\infty}\boldsymbol{p})^{\ast}(\mathscr{L})\in\mathbb{R},
\]
and in the following we will denote by
\[
\int:\overline{\Lambda}{}^{n}(\mathscr{E})\ni\mathscr{L}\longmapsto
\int\mathscr{L}:=[\mathscr{L}]\in\overline{H}{}^{n}(\mathscr{E})
\]
the projection. Thus $\mathscr{L}$ may be interpreted as the Lagrangian
density of a Lagrangian theory constrained by the PDE $\mathscr{E}_{\Phi}$. As
a natural generalization, we interpret $\overline{H}{}(\mathscr{E})$, not only
its $n$-degree component, as space of local function(al)s on $\boldsymbol{M}%
$. By considering all less-dimensional cohomologies rather than just top ones
we have in mind the possibility of defining functionals by integration on
less-dimensional submanifolds of $M$. Such possibility is crucial in
variational calculus with boundary conditions (see \cite{mv07}).

Similarly, for $p>0$, $\mathscr{C}E_{1}^{p,\bullet}(\mathscr{E})$ is naturally
interpreted as space of local differential $p$-forms on $\boldsymbol{M}$.
This informal arguments motivate the

\begin{definition}
Elements in $\overline{H}{}(\mathscr{E})= \mathscr{C}E_{1}^{0,\bullet
}(\mathscr{E})=:\boldsymbol{C}^{\infty}(\boldsymbol{M})^{\bullet}$ are called
\emph{secondary functions on }$\boldsymbol{M}$. For $p>0$, elements in
$H(\mathscr{C}^{p}\Lambda^{p}(\mathscr{E})\otimes\overline{\Lambda}%
{}(\mathscr{E}),\overline{d})\simeq\mathscr{C}E_{1}^{p,\bullet}%
(\mathscr{E})=:\boldsymbol{\Lambda}^{p}(\boldsymbol{M})^{\bullet}$ are called
\emph{secondary differential }$p$\emph{-forms on }$\boldsymbol{M}$. We put
also $\boldsymbol{\Lambda}(\boldsymbol{M})^{\bullet}:=\bigoplus_{p}%
\boldsymbol{\Lambda}^{p}(\boldsymbol{M})^{\bullet}$.
\end{definition}

Notice that elements in $\boldsymbol{\Lambda}(\boldsymbol{M})^{n}$ are
sometimes referred to in the literature as \emph{variational forms} \cite{o86}.

We apply similar arguments to motivate the definition of secondary vector
fields. First of all, notice that there exists a complex
\begin{equation}
\xymatrix@C=17pt{0 \ar[r] & V\mathrm{D}(\mathscr{E}) \ar[r]^-{\overline{S}} & \cdots \ar[r] & 
V\mathrm{D}(\mathscr{E})\otimes\overline{\Lambda}{}^q(\mathscr{E}) \ar[r]^-{\overline{S}} &
V\mathrm{D}(\mathscr{E})\otimes\overline{\Lambda}{}^{q+1}(\mathscr{E}) \ar[r]^-{\overline{S}} & \cdots
}, \label{SpencComp}%
\end{equation}
somehow \textquotedblleft dual" to complex $(\mathscr{C}\Lambda^{1}%
(\mathscr{E})\otimes\overline{\Lambda}{}(\mathscr{E}),\overline{d}%
)\simeq(\mathscr{C}E_{0}^{1,\bullet}(\mathscr{E}),d_{0}^{1,\bullet})$, well
defined by putting
\[
\overline{S}(X\otimes\overline{\omega}):=\overline{S}(X)\wedge\overline
{\omega}+X\otimes\overline{d}\overline{\omega},
\]
$X\in V\mathrm{D}(\mathscr{E})$, $\overline{\omega}\in\overline{\Lambda}%
{}(\mathscr{E})$, where $\overline{S}(X)\in V\mathrm{D}(\mathscr{E})\otimes
\overline{\Lambda}{}^{1}(\mathscr{E})$ is the $V\mathrm{D}(\mathscr{E})$%
-valued horizontal $1$-form defined by putting $\overline{S}%
(X)(Y):=[Y,X]^{V}$, and $[Y,X]^{V}$ is the vertical component of $[Y,X]$.
Complex (\ref{SpencComp}) is called the (horizontal) Spencer complex of
$\mathscr{E}$. As we will see later on in more details, $0$-cohomology
$H^{0}(V\mathrm{D}(\mathscr{E})\otimes\overline{\Lambda}{}%
(\mathscr{E}),\overline{S})$ of the Spencer complex is given by $V\mathrm{D}%
_{\mathscr{C}}(\mathscr{E})$. Now, let $\boldsymbol{X}\in\mathrm{Sym}%
(\mathscr{E})$ and $\rE_{\chi}\in V\mathrm{D}_{\mathscr{C}}(\mathscr{E})$ be
the associated evolutionary vector field, $\chi\in\varkappa|_{\mathscr{E}}$
being a generating section such that $\ell_{\Phi}\chi=0$. Suppose temporarily
that $\rE_{\chi}$ generates a flow $\{A_{t}\}_{t}$ of local diffeomorphisms of
$\mathscr{E}$. Then for any $t$, $A_{t}$ preserves the Cartan distribution and
therefore the image $A_{t}(L)$ of an $n$-dimensional integral submanifold $L$
is an $n$-dimensional integral submanifold. We conclude that $\boldsymbol{X}$
generates a flow of solutions of $\mathscr{E}_{\Phi}$ and, therefore, may be
interpreted as a (local) vector field on $\boldsymbol{M}$. This makes it
rigorous the assertion that \emph{tangent vectors to the solution space of a
PDE are solutions of the associated linearized PDE}. As a natural
generalization, we interpret the whole $H(V\mathrm{D}(\mathscr{E})\otimes
\overline{\Lambda}{}(\mathscr{E}),\overline{S})$, not only its $0$-degree
component, as space of vector fields on $\boldsymbol{M}$. This motivates the

\begin{definition}
Elements in $H(V\mathrm{D}(\mathscr{E})\otimes\overline{\Lambda}%
{}(\mathscr{E}),\overline{S})=:\mathbf{D}(\boldsymbol{M})^{\bullet}$ are
called \emph{secondary vector fields on }$\boldsymbol{M}$.
\end{definition}

All standard operations with vector fields and differential forms have their
secondary analogue. Namely, let $\boldsymbol{\omega}\in\boldsymbol{\Lambda
}^{p}(\boldsymbol{M})^{q}$, $\boldsymbol{\omega}_{1}\in\boldsymbol{\Lambda
}^{p_{1}}(\boldsymbol{M})^{q_{1}}$, $\boldsymbol{\omega}_{2}\in
\boldsymbol{\Lambda}^{p_{2}}(\boldsymbol{M})^{q_{2}}$, $\boldsymbol{X}%
\in\mathbf{D}(\boldsymbol{M})^{r}$, $\boldsymbol{X}_{1}\in\mathbf{D}%
(\boldsymbol{M})^{r_{1}}$, $\boldsymbol{X}_{2}\in\mathbf{D}(\boldsymbol{M}%
)^{r_{2}}$. Then $\boldsymbol{\omega}=[\omega]$, $\boldsymbol{\omega}%
_{1}=[\omega_{1}]$ and $\boldsymbol{\omega}_{2}=[\omega_{2}]$ for some
$\omega\in\mathscr{C}^{p}\Lambda^{p}(\mathscr{E})\otimes\overline{\Lambda}{}%
{}^{q}(\mathscr{E})$, $\omega_{1}\in\mathscr{C}^{p_{1}}\Lambda^{p_{1}%
}(\mathscr{E})\otimes\overline{\Lambda}{}{}^{q_{1}}(\mathscr{E})$ and
$\omega_{2}\in\mathscr{C}^{p_{2}}\Lambda^{p_{2}}(\mathscr{E})\otimes
\overline{\Lambda}{}{}^{q_{2}}(\mathscr{E})$ such that $\overline{d}%
\omega=\overline{d}\omega_{1}=\overline{d}\omega_{2}=0$. Similarly,
$\boldsymbol{X}=[X]$, $\boldsymbol{X}_{1}=[X_{1}]$ and $\boldsymbol{X}%
_{2}=[X_{2}]$ for some $X\in V\mathrm{D}(\mathscr{E})\otimes\overline{\Lambda
}{}{}^{r}(\mathscr{E})$, $X_{1}\in V\mathrm{D}(\mathscr{E})\otimes
\overline{\Lambda}{}{}^{r_{1}}(\mathscr{E})$ and $X_{2}\in V\mathrm{D}%
(\mathscr{E})\otimes\overline{\Lambda}{}{}^{r_{2}}(\mathscr{E})$ such that
$\overline{S}(X)=\overline{S}(X_{1})=\overline{S}(X_{2})=0$. The following
operations are well defined:

exterior product of differential forms:%
\[
\boldsymbol{\omega}_{1}\wedge\boldsymbol{\omega}_{2}:=[(-1)^{q_{1}p_{2}}%
\omega_{1}\wedge\omega_{2}]\in\boldsymbol{\Lambda}^{p_{1}+p_{2}}%
(\boldsymbol{M})^{q_{1}+q_{2}};
\]
exterior differential of a differential form:%
\[
\boldsymbol{d\omega}:=[d^{V}\omega]\in\boldsymbol{\Lambda}^{p+1}%
(\boldsymbol{M})^{q};
\]
commutator of vector fields:%
\[
\lbrack\boldsymbol{X}_{1},\boldsymbol{X}_{2}]:=[[\![X_{1},X_{2}]\!]]\in
\mathbf{D}(\boldsymbol{M})^{r_{1}+r_{2}};
\]
insertion of a vector field into a differential form:%
\[
i_{\boldsymbol{X}}\boldsymbol{\omega}:=[(-1)^{r(p-1)}i_{X}\omega
]\in\boldsymbol{\Lambda}^{p-1}(\boldsymbol{M})^{q+r};
\]
Lie derivative of a differential form along a vector field:%
\[
\mathcal{L}_{\boldsymbol{X}}\boldsymbol{\omega}:=(i_{\boldsymbol{X}}%
\circ\boldsymbol{d}+\boldsymbol{d}\circ i_{\boldsymbol{X}})\boldsymbol{\omega
}\in\boldsymbol{\Lambda}^{p}(\boldsymbol{M})^{q+r};
\]
$[\![{}\cdot{},{}\cdot{}]\!]$ being the Fr\"{o}licher-Nijenhuis bracket of
form-valued vector fields.

Secondary analogue of the standard relations among the above operations hold.
Indeed, let $\boldsymbol{\omega}_{1},\boldsymbol{\omega}_{2},\boldsymbol{X}%
,\boldsymbol{X}_{1},\boldsymbol{X}_{2}$ be as above. The exterior product endows
$\boldsymbol{\Lambda}(\boldsymbol{M})^{\bullet}=\bigoplus_{p,q}%
\boldsymbol{\Lambda}^{p}(\boldsymbol{M})^{q}$ with the structure of a
bi-graded algebra. Namely, $\boldsymbol{\omega}_{1}\wedge\boldsymbol{\omega
}_{2}=(-1)^{p_{1}p_{2}+q_{1}q_{2}}\boldsymbol{\omega}_{2}\wedge
\boldsymbol{\omega}_{1}$. The exterior differential is a bi-graded derivation
of bi-degree $(1,0)$. Namely, $\boldsymbol{d}(\boldsymbol{\omega}_{1}%
\wedge\boldsymbol{\omega}_{2})=\boldsymbol{d\omega}_{1}\wedge
\boldsymbol{\omega}_{2}+(-1)^{p_{1}}\boldsymbol{\omega}_{1}\wedge
\boldsymbol{d\omega}_{2}$. The commutator endows $\mathbf{D}(\boldsymbol{M}%
)^{\bullet}=\bigoplus_{r}$ $\mathbf{D}(\boldsymbol{M})^{r}$ with the structure
of a graded Lie algebra, in particular, $[\boldsymbol{X},[\boldsymbol{X}%
_{1},\boldsymbol{X}_{2}]]=[[\boldsymbol{X},\boldsymbol{X}_{1}],\boldsymbol{X}%
_{2}]+(-1)^{rr_{1}}[\boldsymbol{X}_{1},[\boldsymbol{X},\boldsymbol{X}_{2}]]$.
The `insertion of' and the `Lie derivative along' $\boldsymbol{X}$ are
bi-graded derivations of bi-degree $(-1,r)$ and $(0,r)$ respectively.
Namely, $i_{\boldsymbol{X}}(\boldsymbol{\omega}_{1}\wedge\boldsymbol{\omega
}_{2})=i_{\boldsymbol{X}}\boldsymbol{\omega}_{1}\wedge\boldsymbol{\omega}%
_{2}+(-1)^{p_{1}+rq_{1}}\boldsymbol{\omega}_{1}\wedge i_{\boldsymbol{X}%
}\boldsymbol{\omega}_{2}$ and $\mathcal{L}_{\boldsymbol{X}}(\boldsymbol{\omega
}_{1}\wedge\boldsymbol{\omega}_{2})=\mathcal{L}_{\boldsymbol{X}}%
\boldsymbol{\omega}_{1}\wedge\boldsymbol{\omega}_{2}+(-1)^{rq_{1}%
}\boldsymbol{\omega}_{1}\wedge\mathcal{L}_{\boldsymbol{X}}\boldsymbol{\omega
}_{2}$. Moreover, $[\boldsymbol{d},\boldsymbol{d}]=[\boldsymbol{d}%
,\mathcal{L}_{\boldsymbol{X}}]=[i_{\boldsymbol{X}_{1}},i_{\boldsymbol{X}_{2}%
}]=0$, $[\boldsymbol{d},i_{\boldsymbol{X}}]=\mathcal{L}_{\boldsymbol{X}}$,
$[i_{\boldsymbol{X}_{1}},\mathcal{L}_{\boldsymbol{X}_{2}}]=i_{[\boldsymbol{X}%
_{1},\boldsymbol{X}_{2}]}$, $[\mathcal{L}_{\boldsymbol{X}_{1}},\mathcal{L}%
_{\boldsymbol{X}_{2}}]=\mathcal{L}_{[\boldsymbol{X}_{1},\boldsymbol{X}_{2}]}$,
where $[{}\cdot{},{}\cdot{}]$ denotes the bi-graded commutator.

Despite some time has passed since they were introduced \cite{k92,v84}, to the
author knowledge, no general techniques have been developed so far in order to
effectively compute secondary differential form and vector field spaces, i.e.,
cohomologies of complexes (\ref{HorDRComp}) and (\ref{SpencComp}), in full
generality, other than the one based on the so-called \emph{compatibility
complexes} \cite{t91,v98'} (and, possibly, the Koszul-Tate resolution
\cite{v02}), which is reviewed in the next two sections.

\subsection{Horizontal Calculus on PDEs\label{HorCalc}}

The Cartan distribution determines a \textquotedblleft horizontal differential
calculus\textquotedblright\ on $\mathscr{E}$. Informally speaking, the
horizontal differential calculus is obtained replacing standard partial
derivatives with total derivatives. For instance, a horizontal linear
differential operator is one which is a linear combination of compositions of
total derivatives.

More rigorously, let $\tau:T\longrightarrow\mathscr{E}$ (resp.~$\rho
:R\longrightarrow$ $\mathscr{E}$) be a finite dimensional vector bundle and
$P:=\Gamma(\tau)$ (resp.~$Q:=\Gamma(\rho)$) the $C^{\infty}(\mathscr{E})$%
-module of sections of $\tau$ (resp.~$\rho$). In the following any such
module will be called a \emph{smooth module}. A linear differential operator
$\square:P\longrightarrow Q$ is called a \emph{horizontal (linear)
differential operator} iff, by definition, for any $\theta\in\mathscr{E}$ and
any submanifold $L\subset\mathscr{E}$ such that $\theta\in L$ and $T_{\theta
}L\subset\mathscr{C}_{\theta}$ there exists a differential operator
$\square_{\theta}^{L}:P|_{L}\longrightarrow Q|_{L}$ such that $(\square
p)(\theta)=\square_{\theta}^{L}(p|_{L})(\theta)$ for all $p\in P$. As
examples, notice that horizontal vector fields, the horizontal de Rham
differential $\overline{d}$, the Spencer differential $\overline{S}$ and
universal linearizations are horizontal differential operators. Indeed, Let
$\ldots,e_{a},\ldots$ (resp. $\ldots,\varepsilon_{A},\ldots$) be a local basis
of $P$ (resp.~$Q$). Then a horizontal differential operator $\square
:P\longrightarrow Q$ is characterized as being one locally given by
\begin{equation}
\square p=\varepsilon_{A}\square_{a}^{A}{}^{I}D_{I}^{\mathscr{E}}p^{a}%
,\quad\ldots,\square_{a}^{A}{}^{I},\ldots\ \text{being local functions on
}\mathscr{E}\text{,} \label{LocDiffOp}%
\end{equation}
for all $p=p^{a}e_{a}$ local sections of $\tau$, $\ldots,p^{a},\ldots$ local
functions on $\mathscr{E}$. In particular, if $\mathscr{E}=J^{\infty}$ and
$\mathscr{F}\subset J^{\infty}$ is the infinite prolongation of a PDE, then
any horizontal differential operator $\square:P\longrightarrow Q$ restricts to
$\mathscr{F}$, i.e., there exists a unique (horizontal) differential operator
$\square^{\mathscr{F}}:P|_{\mathscr{F}}\longrightarrow Q|_{\mathscr{F}}$ such
that $\square^{\mathscr{F}}(p|_{\mathscr{F}})=(\square p)|_{\mathscr{F}}$ for
all $p\in P$.

Denote by $\mathscr{C}\mathrm{Diff}(P,Q)$ the set of all horizontal
differential operators $\square:P\longrightarrow Q$. Clearly,
$\mathscr{C}\mathrm{Diff}(P,Q)$ is a $C^{\infty}(\mathscr{E})$-module
naturally isomorphic to $\mathscr{C}\mathrm{Diff}(P,C^{\infty}%
(\mathscr{E}))\otimes Q$ and in what follows we will understand such isomorphism.

Similarly, one may define horizontal jets of sections of vector bundles over
$\mathscr{E}$ just replacing partial derivatives with total derivatives in the
standard definition. We refer to \cite{v02} for the details of the
construction. Analogously to the standard case, one may also define (systems of
horizontal) PDEs determined by linear horizontal differential operators and,
in particular, formally integrable PDEs.

Denote by $\overline{\tau}_{\infty}:\overline{J}{}^{\infty}\tau\longrightarrow
\mathscr{E}$ the bundle of horizontal infinite jets of sections of $\tau$ and
put $\overline{J}{}^{\infty}P:=\Gamma(\overline{\tau}_{\infty})$. For any
$p\in P$ denote by $\overline{j}{}_{\infty}p\in\overline{J}{}^{\infty}P$ its
infinite horizontal jet prolongation. There is a canonical monomorphism of
$C^{\infty}(\mathscr{E})$-modules $h:\mathscr{C}\mathrm{Diff}(P,Q)\ni
\square\longmapsto h_{\square}\in\mathrm{Hom}(\overline{J}{}^{\infty}P,Q)$,
where $h_{\square}$ is the unique $C^{\infty}(\mathscr{E})$-linear map such
that $h_{\square}(\overline{j}{}_{\infty}p)=\square p$ for all $p\in P$.
Moreover $h_{\square}$ can be uniquely prolonged to a $C^{\infty
}(\mathscr{E})$-linear map $h_{\square}^{\infty}:\overline{J}{}^{\infty
}P\longrightarrow\overline{J}{}^{\infty}Q$ such that $h_{\square}^{\infty
}(\overline{j}{}_{\infty}p)=\overline{j}{}_{\infty}(\square p)$ for all $p\in
P$.

The following remarkable correspondence,
\begin{equation}
\overline{J}{}^{\infty}\varkappa\ni\overline{j}{}_{\infty}\chi\longmapsto
\rE_{\chi}\in V\mathrm{D}(J^{\infty}), \label{IsoEtaJ}%
\end{equation}
determines a well defined isomorphism of $C^{\infty}(J^{\infty})$-modules.
The dual isomorphism is given by
\begin{equation}
\mathscr{C}\Lambda^{1}(J^{\infty})\ni\omega\longmapsto\square_{\omega}%
\in\mathscr{C}\mathrm{Diff}(\varkappa,C^{\infty}(J^{\infty})),
\label{IsoEtaDualJ}%
\end{equation}
where $\square_{\omega}:\varkappa\longrightarrow C^{\infty}(J^{\infty})$ is
defined by putting $\square_{\omega}\chi:=\omega(\rE_{\chi})$, $\chi
\in\varkappa$. Accordingly, there is a natural embedding $\eta_{\Phi
}:V\mathrm{D}(\mathscr{E})\hookrightarrow\overline{J}{}^{\infty}%
\varkappa|_{\mathscr{E}}$ given by the composition
\[
\xymatrix{V\mathrm{D}(\mathscr{E})\,\ar@{^{(}->}[r] \ar@/_1.2pc/[rr]_-{\eta_\Phi} & V\mathrm{D}(J^\infty)|_{\mathscr{E}} \ar@{}@<-0.8ex>[r]^-{\widetilde{\quad\quad}} \ar[r] & \overline{J}{}^\infty \varkappa |_{\mathscr{E}} },
\]
and, dually, a natural projection $\eta_{\Phi}^{\ast}:\mathscr{C}\mathrm{Diff}%
(\varkappa|_{\mathscr{E}},C^{\infty}(\mathscr{E}))\twoheadrightarrow
\mathscr{C}\Lambda^{1}(\mathscr{E})$ given by the composition
\[
\xymatrix{\mathscr{C}\mathrm{Diff}(\varkappa|_\mathscr{E},C^\infty(\mathscr{E}))\ar[r] \ar@{}@<-0.8ex>[r]^-{\widetilde{\quad\quad}} \ar[r] \ar@/_1.3pc/[rr]_-{\eta^\ast_\Phi} & \mathscr{C}\Lambda^1(J^\infty)|_{\mathscr{E}}  \ar@{->>}[r] & \mathscr{C}\Lambda^1(\mathscr{E}) },
\]
where the arrows \textquotedblleft$\widetilde{\longrightarrow}$%
\textquotedblright\ are the inverses of restrictions to $\mathscr{E}$ of
isomorphisms (\ref{IsoEtaJ}) and (\ref{IsoEtaDualJ}), respectively. Finally, let $P$ be the smooth module where $\Phi$ belongs. Notice that
the sequence
\begin{equation}
\xymatrix{0 \ar[r] & V\mathrm{D}(\mathscr{E}) \ar[r]^-{\eta_\Phi} & \overline{J}{}^\infty \varkappa |_\mathscr{E} \ar[r]^-{h_\Phi^\infty} & \overline{J}{}^\infty P }, \label{ExSeq}%
\end{equation}
where $h_{\Phi}:=h_{\ell_{\Phi}}$, and its dual
\begin{equation}
\xymatrix{ \mathscr{C}\mathrm{Diff}(P,C^\infty(\mathscr{E}) )\ar[r]^-{h_\Phi^\infty{}^\ast} & \mathscr{C}\mathrm{Diff}(\varkappa |_\mathscr{E},C^\infty(\mathscr{E}) ) \ar[r]^-{\eta_\Phi^\ast}  &  \mathscr{C}\Lambda^1(\mathscr{E}) \ar[r] & 0}, \label{ExSeqDual}
\end{equation}
where $h_{\Phi}^{\infty}{}^{\ast}(\Delta):=\Delta\circ\ell_{\Phi}$, $\Delta
\in\mathscr{C}\mathrm{Diff}(P,C^{\infty}(\mathscr{E}))$, are exact.

There exists a horizontal analogue of the concept of adjoint operator to a
linear differential operator. Let $R$ be a smooth module (see above). Put
$R^{\dag}:=\mathrm{Hom}(R,\overline{\Lambda}{}
^{n}(\mathscr{E}))$. $R^{\dag}$ is a smooth module as well and it is called
the \emph{adjoint module to }$R$. Obviously, $R^{\dag\dag}$ identifies
canonically with $R$. Denote by $R^{\dag}\times R\ni(r^{\dag},r)\longmapsto
\langle r^{\dag},r\rangle:=r^{\dag}(r)\in\overline{\Lambda}{}^{n}%
(\mathscr{E})$ the natural bi-linear pairing. For any local basis
$\ldots,\kappa_{a},\ldots$ of $R$ we denote by $\ldots,\kappa^{\dag}{}%
^{a},\ldots$ the local basis of $R^{\dag}$ such that $\kappa^{\dag}{}^{a}$ is
the local homomorphism $R\longrightarrow\overline{\Lambda}{}^{n}(\mathscr{E})$
defined by putting $\langle\kappa^{\dag}{}^{a},\kappa_{b}\rangle:=\delta
_{b}^{a}\overline{\mathrm{d}}{}^{n}x$ and $\overline{\mathrm{d}}{}%
^{n}x:=\overline{d}x^{1}\wedge\cdots\wedge\overline{d}x^{n}$, $a,b=1,2,\ldots$.

\begin{proposition}
\label{DuBois-Reymond} Let $r\in R$ (resp.~$r^{\dag}\in R^{\dag}$), then $r=0$
(resp.~$r^{\dag}=0$) iff $\int\langle r^{\dag},r\rangle=0$ for all $r^{\dag
}\in R^{\dag}$ (resp.~$r\in R$).
\end{proposition}

Proposition \ref{DuBois-Reymond} may be referred to as the \emph{cohomological
DuBois-Reymond theorem }and will be used later on without further comments.

Now let $P,Q$ be as above and $\square:P\longrightarrow Q$ a horizontal
differential operator. It can be proved that there exists a unique
differential operator (of the same order as $\square$) $\square^{\dag}%
:Q^{\dag}\longrightarrow P^{\dag}$ such that%
\begin{equation}
\int\langle q^{\dag},\square p\rangle=\int\langle\square^{\dag}q^{\dag
},p\rangle\label{Green}%
\end{equation}
for all $p\in P$, $q^{\dag}\in Q^{\dag}$. $\square^{\dag}$ is called the
\emph{adjoint operator to }$\square$ and (\ref{Green}) is called the
(horizontal) \emph{Green formula }\cite{b...99,k92,v01}.

Adjoint operators have the following properties. First, $\square^{\dag\dag
}=\square$. Second, let $\Delta:Q\longrightarrow R$ be another horizontal differential
operator, then $(\Delta\circ\square)^{\dag}=\square^{\dag}\circ\Delta^{\dag}$.
If $\square$ is locally given by (\ref{LocDiffOp}) then $\square^{\dag}$ is
locally given by
\[
\square^{\dag}q^{\dag}=(-1)^{|I|}e^{\dag}{}^{a}D_{I}(\square_{a}^{A}{}%
^{I}q_{A}^{\dag}),
\]
for all $q^{\dag}=q_{A}^{\dag}\varepsilon^{\dag}{}^{A}$ local elements of
$Q^{\dag}$, $\ldots,q_{A}^{\dag},\ldots$ local functions on $\mathscr{E}$.
As an example, notice that the adjoint module of $\overline{\Lambda}{}%
^{q}(\mathscr{E})$ is canonically isomorphic to $\overline{\Lambda}{}%
^{n-q}(\mathscr{E})$, and that the adjoint operator of the horizontal de Rham
differential $\overline{d}:\overline{\Lambda}{}^{q}%
(\mathscr{E})\longrightarrow\overline{\Lambda}{}^{q+1}(\mathscr{E})$ is the
operator $(-1)^{n-q-1}\overline{d}:\overline{\Lambda}{}^{n-q-1}%
(\mathscr{E})\longrightarrow\overline{\Lambda}{}^{n-q}(\mathscr{E})$,
$q=0,\ldots,n$.

Notice that the Green formula amounts to say that for any $p\in P$, $q^{\dag
}\in Q^{\dag}$ there exists $\rG_{p,q^{\dag}}\in\overline{\Lambda}{}%
^{n-1}(\mathscr{E})$ such that $\langle q^{\dag},\square p\rangle
-\langle\square^{\dag}q^{\dag},p\rangle=\overline{d}\rG_{p,q^{\dag}}$. It can
be proved \cite{a02} that $\rG_{p,q^{\dag}}$ can be chosen of the form
$\rG(p,q^{\dag})$, $\rG:P\times Q^{\dag}\longrightarrow\overline{\Lambda}%
{}^{n-1}(\mathscr{E})$ being a (possibly non unique) horizontal
bi-differential operator independent of $p$ and $q^{\dag}$. Any such operator
$\rG$ is called a \emph{Legendre operator for }$\square$ \cite{av04}. The
Green formula plays a central role in the theory of the $\mathscr{C}$%
-spectral sequence.

\subsection{Formal Theory of Horizontal PDEs and Secondary
Calculus\label{SpencGold}}

There exists a horizontal analogue of the Goldschmidt-Spencer formal theory
of linear differential equations (see \cite{g67A,s69} for a complete account
of the classical theory - see also \cite{g67B} - and \cite{kv98,v98'} for
its horizontal analogue).

Let $\Delta:P\longrightarrow P_{1}$ be a horizontal differential operator of
order $\leq k$ between smooth modules.

\begin{definition}
\label{Comp}A complex of horizontal differential operators between smooth
modules
\begin{equation}
\xymatrix@C=30pt{0 \ar[r] & P \ar[r]^-\Delta & P_1 \ar[r]^-{\Delta_1} & \cdots \ar[r] & P_q \ar[r]^-{\Delta_q}  & P_{q+1} \ar[r]^-{\Delta_{q+1}} & \cdots} \label{CompComp}%
\end{equation}
is called a \emph{compatibility complex for} $\Delta$ iff the sequence of
homomorphisms%
\[
\xymatrix@C=35pt{0 \ar[r] & \overline{J}{}^\infty P \ar[r]^-{h_\Delta^\infty }& \overline{J}{}^\infty P_1 \ar[r]^-{h_{\Delta_1}^\infty } & \cdots \ar[r] & \overline{J}{}^\infty P_q \ar[r]^-{h_{\Delta_q}^\infty }  & \overline{J}{}^\infty P_{q+1} \ar[r]^-{h_{\Delta_{q+1}}^\infty } & \cdots}
\]
is exact. $\Delta_{1}$ is called a \emph{compatibility operator for} $\Delta$.
\end{definition}

The existence of a non trivial compatibility operator for $\Delta$ formalizes
the fact that the equation $\Delta p=0$ is overdetermined \cite{s69}. We
stress that Definition \ref{Comp} is slightly different from the one usually
found in the literature (see, for instance, \cite{g67A,kv98}). However, it can
be shown that, if $\Delta$ determines a formally integrable PDE, then the two
coincide, and Definition \ref{Comp} is the most suitable for our purposes.

\begin{theorem}
[Goldschmidt]Let $\Delta$ be a horizontal differential operator between smooth
modules. If $\Delta$ determines a formally integrable horizontal PDE, then
there exists a (non unique) compatibility complex (\ref{CompComp}) for $\Delta
$, such that $\Delta_{i}$ determines a formally integrable horizontal PDE for
any $i=1,2,\ldots$.
\end{theorem}

Any compatibility complex as in the above theorem will be said \emph{regular}.
Let $\Delta:P\longrightarrow P_{1}$ determine a formally integrable PDE and
(\ref{CompComp}) be a regular compatibility complex for it. Then, the
compatibility operator $\Delta_{1}$ has the following remarkable property.

\begin{proposition}
\label{Prop1}Let $\square:P_{1}\longrightarrow Q$ be a horizontal differential
operator such that $\square\circ\Delta=0$. Then there exists a horizontal
differential operator $\nabla:P_{2}\longrightarrow Q$ such that $\square
=\nabla\circ\Delta_{1}$. If $\Delta_{2}=0$ then $\nabla$ is unique.
\end{proposition}

Thus, let $\Delta:P\longrightarrow P_{1}$ determine a formally integrable PDE
and
\[
\xymatrix@C=30pt{0 \ar[r] & P \ar[r]^-\Delta & P_1 \ar[r]^-{\Delta_1} & \cdots \ar[r] & P_{s-1} \ar[r]^-{\Delta_{s-1}}  & P_{s} \ar[r] & 0}
\]
be a finite length regular compatibility complex. In this situation we say
that \emph{the compatibility length of }$\Delta$\emph{ is} $\leq s$.

Now let $\pi:E\longrightarrow M$ be a fiber bundle, $\tau:T\longrightarrow
J^{\infty}$ a vector bundle, $\Phi\in\mathrm{diff}(\pi,\tau)$ and
$\mathscr{E}:=\mathscr{E}_{\Phi}^{(\infty)}$. Put $P_{1}:=\Gamma
(\tau)|_{\mathscr{E}}$. Notice that if $\mathscr{E}_{\Phi}$ is a formally
integrable PDE, then $\ell_{\Phi}:\varkappa|_{\mathscr{E}}\longrightarrow
P_{1}$ determines a formally integrable, linear, horizontal PDE \cite{t91}.

\begin{theorem}
[Spencer]\label{SpencTh}Cohomology $\mathbf{D}(\boldsymbol{M})^{\bullet}$ of
complex $(V\mathrm{D}(\mathscr{E})\otimes\overline{\Lambda}%
(\mathscr{E}),\overline{S})$ is canonically isomorphic to cohomology of any
regular compatibility complex
\[
\xymatrix@C=30pt{0 \ar[r] & \varkappa |_{\mathscr{E}} \ar[r]^-{\ell_\Phi} & P_1 \ar[r]^-{\Delta_1} & \cdots \ar[r] & P_q \ar[r]^-{\Delta_q}  & P_{q+1} \ar[r]^-{\Delta_{q+1}} & \cdots}
\]
for $\ell_{\Phi}$.
\end{theorem}

In the following we will only consider regular compatibility complexes.

Isomorphism $\ker\ell_{\Phi}\simeq\mathbf{D}(\boldsymbol{M})^{0}$ is given by
\[
\ker\ell_{\Phi}\ni\chi\longmapsto\rE_{\chi}\in V\mathrm{D}_{\mathscr{C}}%
(\mathscr{E})=\mathbf{D}(\boldsymbol{M})^{0}.
\]

We now describe isomorphism $\ker\Delta_{1}/\operatorname{im}\ell_{\Phi}%
\simeq\mathbf{D}(\boldsymbol{M})^{1}$ referring to \cite{kv98} for the
remaining homogeneous components. Let $p\in P_{1}$ be such that $\Delta
_{1}p=0$. Consider $\overline{j}{}_{\infty}p\in\overline{J}{}^{\infty}P_{1}$.
Then $h_{\Delta_{1}}(\overline{j}{}_{\infty}p)=\Delta_{1}p=0$ and, therefore,
$\Delta_{1}$ being a compatibility operator for $\ell_{\Phi}$, there exists
$j\in\overline{J}{}^{\infty}\varkappa$ such that $\overline{j}{}_{\infty
}p=h_{\Phi}^\infty(j|_{\mathscr{E}})$. Let $X:=\eta_{0}^{-1}(j)\in V\mathrm{D}%
(J^{\infty})$ (here $0$ is the trivial differential operator) and
$\widetilde{\Omega}:=\overline{S}(X)\in V\mathrm{D}(J^{\infty})\otimes
\overline{\Lambda}{}^{1}(J^{\infty})$. It is easy to prove, suitably using
exactness of sequence (\ref{ExSeq}), that $\widetilde{\Omega}$ restricts to
$\mathscr{E}$, i.e., $\Omega:=\widetilde{\Omega}|_{\mathscr{E}}\in
V\mathrm{D}(\mathscr{E})\otimes\overline{\Lambda}{}^{1}(\mathscr{E})\subset
V\mathrm{D}(J^{\infty})|_{\mathscr{E}}\otimes\overline{\Lambda}{}%
^{1}(\mathscr{E})$. Moreover, $\overline{S}(\Omega)=0$. Finally, the
isomorphism $\ker\Delta_{1}/\operatorname{im}\ell_{\Phi}\simeq\mathbf{D}%
(\boldsymbol{M})^{1}$ maps $p+\operatorname{im}\ell_{\Phi}^{\mathscr{E}}%
\in\ker\Delta_{1}/\operatorname{im}\ell_{\Phi}$ to $[\Omega]\in H^{1}%
(V\mathrm{D}(\mathscr{E})\otimes\overline{\Lambda}{}(\mathscr{E}),\overline
{S})=\mathbf{D}(\boldsymbol{M})^{1}$.

\begin{corollary}
Cohomology $\mathbf{\Lambda}^{1}(\boldsymbol{M})^{\bullet}$ of complex
$(\mathscr{C}\Lambda^{1}(\mathscr{E})\otimes\overline{\Lambda}%
(\mathscr{E}),\overline{d})$ is canonically isomorphic to homology of the
adjoint complex
\[
\xymatrix@C=30pt{0  & \varkappa |_{\mathscr{E}}^\dag \ar[l]  & P_1^\dag \ar[l]_-{\ell_\Phi^{\dag}} & \cdots \ar[l] & P_q^\dag \ar[l]_-{\Delta_{q-1}^\dag}  & P_{q+1}^\dag  \ar[l]_-{\Delta_{q}^\dag}& \cdots \ar[l]}
\]
of any regular compatibility complex for $\ell_{\Phi}$.
\end{corollary}

Isomorphism $\mathbf{\Lambda}^{1}(\boldsymbol{M})^{n}\simeq
\operatorname{coker}\ell_{\Phi}^{\dag}$ is described as follows. Projection
$\eta_{\Phi}^{\ast}:\mathscr{C}\mathrm{Diff}(\varkappa|_{\mathscr{E}}%
,C^{\infty}(\mathscr{E}))\twoheadrightarrow\mathscr{C}\Lambda^{1}%
(\mathscr{E})$ gives rise to a projection
\[
\eta_{\Phi}^{\ast}\otimes\mathrm{id}_{\overline{\Lambda}(\mathscr{E})}%
:\mathscr{C}\mathrm{Diff}(\varkappa|_{\mathscr{E}},\overline{\Lambda}%
{}(\mathscr{E}))\twoheadrightarrow\mathscr{C}\Lambda^{1}(\mathscr{E})\otimes
\overline{\Lambda}{}(\mathscr{E})
\]
which, abusing the notation, we denote again by $\eta_{\Phi}^{\ast}$. Thus, let
$\omega\in\mathscr{C}\Lambda^{1}(\mathscr{E})\otimes\overline{\Lambda}{}%
^{n}(\mathscr{E})$ and $\square\in$ $\mathscr{C}\mathrm{Diff}(\varkappa
|_{\mathscr{E}},\overline{\Lambda}{}^{n}(\mathscr{E}))$ be such that
$\eta_{\Phi}^{\ast}{}(\square)=\omega$. Consider $\square^{\dag}:C^{\infty
}(\mathscr{E})\longrightarrow\varkappa|_{\mathscr{E}}^{\dag}$. Isomorphism
$\mathbf{\Lambda}^{1}(\boldsymbol{M})^{n}\simeq\operatorname{coker}\ell_{\Phi
}^{\dag}$ maps $[\omega]\in H^{n}(\mathscr{C}\Lambda^{1}(\mathscr{E})\otimes
\overline{\Lambda}{}(\mathscr{E}),\overline{d})=\mathbf{\Lambda}%
^{1}(\boldsymbol{M})^{n}$ to $\square^{\dag}1+\operatorname{im}\ell_{\Phi
}^{\dag}\in\operatorname{coker}\ell_{\Phi}^{\dag}{}$.

We now describe isomorphism $\mathbf{\Lambda}^{1}(\boldsymbol{M})^{n-1}%
\simeq\ker\ell_{\Phi}^{\dag}/\operatorname{im}\Delta_{1}^{\dag}$ referring
again to \cite{kv98} for the remaining homogeneous components. Let $\omega
\in\mathscr{C}\Lambda^{1}(\mathscr{E})\otimes\overline{\Lambda}{}%
^{n-1}(\mathscr{E})$ and $\square\in$ $\mathscr{C}\mathrm{Diff}(\varkappa
|_{\mathscr{E}},\overline{\Lambda}{}^{n}(\mathscr{E}))$ be such that
$\overline{d}\omega=0$ and $\eta_{\Phi}^{\ast}{}(\square)=\omega$. Then, it
follows from exactness of sequence (\ref{ExSeqDual}) that $\overline{d}%
\circ\square=\Delta\circ\ell_{\Phi}$ for some $\Delta\in
\mathscr{C}\mathrm{Diff}(P_{1},\overline{\Lambda}{}^{n}(\mathscr{E}))$.
Consider $\Delta^{\dag}:C^{\infty}(\mathscr{E})\longrightarrow P_{1}^{\dag}$
and put $p^{\dag}:=\Delta^{\dag}1\in P_{1}^{\dag}$. We have $\ell_{\Phi}%
^{\dag}(p^{\dag})=(\ell_{\Phi}^{\dag}\circ\Delta^{\dag})(1)=(\Delta\circ
\ell_{\Phi})^{\dag}(1)=(\overline{d}\circ\square)^{\dag}(1)=(\square^{\dag
}\circ\overline{d}{}^{\dag})(1)=(\square^{\dag}\circ d)(1)=0$. Thus $p^{\dag
}\in\ker\ell_{\Phi}^{\dag}$. Isomorphism $\mathbf{\Lambda}^{1}(\boldsymbol{M}%
)^{n-1}\simeq\ker\ell_{\Phi}^{\dag}/\operatorname{im}\Delta_{1}^{\dag}$ maps
$[\omega]\in H^{n-1}(\mathscr{C}\Lambda^{1}(\mathscr{E})\otimes\overline
{\Lambda}{}(\mathscr{E}),\overline{d})=\mathbf{\Lambda}^{1}(\boldsymbol{M}%
)^{n-1}$ to $p^{\dag}+\operatorname{im}\Delta_{1}^{\dag}\in\ker\ell_{\Phi
}^{\dag}/\operatorname{im}\Delta_{1}^{\dag}$.

Notice that the above corollary describes to some extent the $1$-st column of
the $1$-st term of the $\mathscr{C}$-spectral sequence of $\mathscr{E}$. The
following theorem due to Verbovetsky \cite{v98'} (see also \cite{v84} for the case $s=2$) extends
it to the remaining columns.

\begin{theorem}
[$s$-lines]\label{slines}Let $\mathscr{E}\subset J^{\infty}$ be the infinite
prolongation of a formally integrable PDE $\mathscr{E}_{\Phi}$ and let the
compatibility length of $\ell_{\Phi}$ be $\leq s$. Then $\mathscr{C}E_{1}%
^{p,q}(\mathscr{E})=0$ if $p>0$ and $q<n-s$.
\end{theorem}

\begin{example}
[empty equation]\label{Ex1}If $\Phi=0$ then $\mathscr{E}=J^{\infty}$,
$\ell_{\Phi}=0$ and its compatibility length is $0$. In this case
$\mathbf{D}(\boldsymbol{M})^{r}=0$ for $r\neq0$ and $\mathbf{D}(\boldsymbol{M}%
)^{\bullet}=\mathbf{D}(\boldsymbol{M})^{0}\simeq\varkappa$. The exact
sequence
\begin{equation}
\xymatrix{0 \ar[r] & \mathbf{D}(\boldsymbol{M})^0 \ar[r] & V\mathrm{D}(J^\infty) \ar@/^1.1pc/[l]^-{\psi} \ar[r]^-{\overline{S}}& V\mathrm{D}(J^\infty)\otimes\overline{\Lambda}{}^1(J^\infty)} \label{SplitD}%
\end{equation}
splits via the composition
\begin{equation}
\xymatrix{V\mathrm{D}(J^{\infty}) \ar[r]^-{\eta_{0}} \ar@/_1.3pc/[rrr]_-{\psi}& \overline{J}{}^{\infty}(\varkappa) \ar@{->>}[r]&
\varkappa \ar@{}@<-1.0ex>[r]^-{\widetilde{\quad\quad\quad}} \ar@{-}@<-0.1ex>[r]& \mathbf{D}(\boldsymbol{M})^{0}}. \label{SplitL}%
\end{equation}
Similarly, $\mathbf{\Lambda}^{1}(\boldsymbol{M})^{q}=0$ for $q\neq n$ and
$\mathbf{\Lambda}^{1}(\boldsymbol{M})^{\bullet}=\mathbf{\Lambda}%
^{1}(\boldsymbol{M})^{n}\simeq\varkappa^{\dag}$. The exact sequence%
\[
\xymatrix{\mathscr{C}\Lambda^1(J^\infty)\otimes \overline{\Lambda}{}^{n-1}(J^\infty) \ar[r]^-{\overline{d}} & \mathscr{C}\Lambda^1(J^\infty)\otimes \overline{\Lambda}{}^n(J^\infty) \ar[r] & \boldsymbol{\Lambda}^1(\boldsymbol{M})^n \ar@/^1.3pc/[l]^-{\psi^\dag} \ar[r] & 0 }
\]
splits via the composition
\[
\xymatrix{\boldsymbol{\Lambda}^1(\boldsymbol{M})^n  \ar@{}@<-1.0ex>[r]^-{\widetilde{\quad\quad\quad}} \ar@{-}@<-0.2ex>[r] \ar@/_1.7pc/[rrrr]_-{\psi^\dag}& \text{\raisebox{1.5mm}{$\varkappa^\dag$}}\, \ar@<-0.2ex>@{^(->}[r]& \mathscr{C}\mathrm{Diff}(\varkappa,\overline{\Lambda}{}^n) \ar[rr]^-{\eta_{0}^\ast} & & \mathscr{C}\Lambda^1(J^\infty)\otimes\overline{\Lambda}{}^n }.
\]

Let $Y\in V\mathrm{D}(J^{\infty})$ and $\varphi\in\varkappa^{\dag}$ be locally
given by $Y=Y_{I}^{\alpha}\partial_{\alpha}^{I}$ and $\varphi=\varphi_{\alpha
}\partial^{\dag\alpha}$, $\ldots,Y_{I}^{\alpha},\ldots,\varphi_{\alpha}%
,\ldots$ being local functions on $J^{\infty}$. Then, locally,
\[
\psi(Y)=Y_{\mathsf{O}}^{\alpha}\partial_{\alpha}\quad\text{and}\quad
\psi^{\dag}(\varphi)=\varphi_{\alpha}\omega_{\mathsf{O}}^{\alpha}%
\otimes\overline{\mathrm{d}}{}^{n}x.
\]

Finally, notice that both diagrams (\ref{SplitD}) and (\ref{SplitL}) restrict
to the infinite prolongation of a PDE and such restrictions preserve the exactness.
\end{example}

In the following we will understand the above isomorphisms $\mathbf{D}(\boldsymbol{M}%
)^{0}\simeq\varkappa$ and $\mathbf{\Lambda}^{1}(\boldsymbol{M})^{n}%
\simeq\varkappa^{\dag}$. In order not to make the notation too heavy we will
also understand the monomorphism $\psi^{\dag}$. According to this convention
$\varkappa^{\dag}$ is understood as a subset in $\mathscr{C}\Lambda
^{1}(J^{\infty})\otimes\overline{\Lambda}{}^{n}(J^{\infty})$. Moreover, if
$\varphi\in\varkappa^{\dag}$ and $\chi\in\varkappa$, then $i_{\rE_{\chi}}%
\psi^{\dag}(\varphi)\in\overline{\Lambda}{}^{n}(J^{\infty})$ identifies with
$\langle\varphi,\chi\rangle$.

\begin{example}
[irreducible equations]\label{IrrEq}A non-empty PDE $\mathscr{E}_{\Phi}$ is
called $\ell$\emph{-normal} (or, in physical terms, \emph{irreducible}) iff
the compatibility length of $\ell_{\Phi}$ is $\leq1$. In this case $\Delta
_{1}$ may be chosen equal to $0$, $\mathbf{D}(\boldsymbol{M})^{r}=0$ for
$r\neq0,1$, $\mathbf{D}(\boldsymbol{M})^{0}\simeq\ker\ell_{\Phi}$ as above
and
\[
\mathbf{D}(\boldsymbol{M})^{1}\simeq\operatorname{coker}\ell_{\Phi}.
\]
Similarly, $\mathbf{\Lambda}^{1}(\boldsymbol{M})^{q}=0$ for $q\neq n,n-1$,
$\mathbf{\Lambda}^{1}(\boldsymbol{M})^{n}\simeq\operatorname{coker}\ell_{\Phi
}^{\dag}$ as above and
\[
\mathbf{\Lambda}^{1}(\boldsymbol{M})^{n-1}\simeq\ker\ell_{\Phi}^{\dag}.
\]

\end{example}

\section{The Covariant Phase Space}

\subsection{Lagrangian Field Theories and the CPS\label{Sec7}}

The calculus of variations is formalized in a coordinate-free way via the
$\mathscr{C}$-spectral sequence.

\begin{definition}
A \emph{Lagrangian (field) theory} is the datum $(\pi,\boldsymbol{S})$ of a
fiber bundle $\pi:E\longrightarrow M$ and an \emph{action} $\boldsymbol{S}%
\in\overline{H}{}^{n}(J^{\infty})$. Any representative $\mathscr{L}\in
\overline{\Lambda}{}^{n}(J^{\infty})$ of the cohomology class $\boldsymbol{S}%
=\int\mathscr{L}$ is called a \emph{Lagrangian density} of the theory
$(\pi,\boldsymbol{S})$.
\end{definition}

Recall that the space $\boldsymbol{M}$ of $n$-dimensional integral
submanifolds of the Cartan distribution $\mathscr{C}$ on $J^{\infty}$ is in
one-to-one correspondence (via infinite jet prolongation) with the space of
local sections of $\pi$. As above we will often identify $\boldsymbol{M}$ with
the \textquotedblleft store\textquotedblright\ $(J^{\infty},\mathscr{C})$ of
its elements. $\boldsymbol{M}$ is known in the Physics literature as the
\emph{space of histories} and an action $\boldsymbol{S}\in\overline{H}{}%
^{n}(J^{\infty})\subset\boldsymbol{C}^{\infty}(\boldsymbol{M})^{\bullet}$ is a
secondary function on it.

Within secondary calculus, the Euler-Lagrange equations (whose solutions make
it stationary the action) associated to the Lagrangian theory $(\pi
,\boldsymbol{S})$ are easily obtained by applying to $\boldsymbol{S}$ the
secondary de Rham differential $\boldsymbol{d}:\boldsymbol{C}^{\infty
}(\boldsymbol{M})^{\bullet}\longrightarrow\boldsymbol{\Lambda}^{1}%
(\boldsymbol{M})^{\bullet}$. Indeed, according to the previous section,
$\boldsymbol{\Lambda}^{1}(\boldsymbol{M})^{\bullet}=\boldsymbol{\Lambda}%
^{1}(\boldsymbol{M})^{n}\simeq\varkappa^\dag$ and $\boldsymbol{dS}$
identifies with the element $\boldsymbol{E}(\mathscr{L}):=\tilde{\ell}
{}_{\mathscr{L}}^{\,\dag}1\in\varkappa^{\dag}\subset\mathscr{C}\Lambda
^{1}(J^{\infty})\otimes\overline{\Lambda}{}^{n}(J^{\infty})$, where we put
$\tilde{\ell}{}_{\mathscr{L}}:=(\eta_{0}^{\ast})^{-1}(d^{V}%
\mathscr{L}):\varkappa\longrightarrow\overline{\Lambda}{}^{n}(J^{\infty})$,
$\mathscr{L}$ being any Lagrangian density. Locally, $\mathscr{L}=L\,\overline
{\mathrm{d}}{}^{n}x$ for some local function $L=L(\ldots,x^{i},\ldots
,u_{I}^{\alpha},\ldots)$ on $J^{\infty}$ and
\[
\boldsymbol{E}(\mathscr{L})=\tfrac{\delta L}{\delta u^{\alpha}}\partial^{\dag
}{}^{\alpha},
\]
$\tfrac{\delta L}{\delta u^{\alpha}}:=(-1)^{|I|}D_{I}(\partial_{\alpha}^{I}L)$
being the so-called \emph{Euler-Lagrange derivatives} of $L$, $\alpha
=1,\ldots,m$. Thus $\boldsymbol{dS}$ is naturally interpreted as the left hand
side of the Euler-Lagrange equations $\mathscr{E}_{\boldsymbol{E}%
(\mathscr{L})}$ of the theory $(\pi,\boldsymbol{S})$. In the following we will
always assume $\mathscr{E}_{\boldsymbol{E}(\mathscr{L})}$ to be a formally
integrable PDE.

Let $\mathscr{E}:=\mathscr{E}_{\boldsymbol{E}(\mathscr{L})}^{(\infty)}$. The
space $\boldsymbol{P}$ of $n$-dimensional integral submanifolds of the Cartan
distribution on $\mathscr{E}$ is in one-to-one correspondence with the space
of (local) solutions of $\mathscr{E}_{\boldsymbol{E}(\mathscr{L})}$ and is
called, according to Physics literature, the (non-reduced) \emph{CPS} of the
theory $(\pi,\boldsymbol{S})$ \cite{bhs91,abr91,c88,cw87,lw90}.

By definition%
\begin{equation}
d^{V}\mathscr{L}-\boldsymbol{E}(\mathscr{L})=\overline{d}\theta
\label{Legendre}%
\end{equation}
for some $\theta\in\mathscr{C}\Lambda^{1}(J^{\infty})\otimes\overline{\Lambda
}{}^{n-1}(J^{\infty})$. Any such $\theta$ will be called a \emph{Legendre
form} \cite{av04} (notice that $\mathscr{L}-\theta$ is a so-called
\emph{Lepagean equivalent} \cite{k83,k73} of $\mathscr{L}$). Equation (\ref{Legendre}) may be interpreted as the \emph{first variation formula} for the action $\boldsymbol{S}$. In this respect, the existence of a global Legendre form was first discussed in \cite{k80}. Any two Legendre
forms $\theta.\theta^{\prime}$ differ by a closed, and therefore exact, form
$\overline{d}\lambda$, $\lambda\in\mathscr{C}\Lambda^{1}(J^{\infty}%
)\otimes\overline{\Lambda}{}^{n-2}(J^{\infty})$ (see, for instance,
\cite{av04,lw90} for a local description of Legendre forms). Notice that, in
view of isomorphism $\eta_{0}^{\ast}$, identity (\ref{Legendre}) may be
understood as the Green formula
\[
\tilde{\ell}{}_{\mathscr{L}}-\tilde{\ell}{}_{\mathscr{L}}^{\,\dag}%
1=(\overline{d}\circ\rG)({}\cdot{},1)
\]
for the horizontal operator $\tilde{\ell}{}_{\mathscr{L}}$, $\rG$ being a
Legendre operator for it.

\begin{theorem}
[Zuckerman]There is a closed secondary $2$-form $\boldsymbol{\omega}$ on
$\boldsymbol{P}$ canonically determined by the corresponding Lagrangian theory
$(\pi,\boldsymbol{S})$.
\end{theorem}

\begin{proof}
Using the language introduced so far we reproduce here the proof in \cite{z87}
by adding the only missing point, that is the independence of
$\boldsymbol{\omega}$ of the choice of a Lagrangian density. Thus, let
$\theta$ be a Legendre form. Put
\[
\omega:=-i_{\mathscr{E}}^{\ast}(d^{V}\theta)\in\mathscr{C}^2\Lambda
^{2}(\mathscr{E})\otimes\overline{\Lambda}{}^{n-1}(\mathscr{E}).
\]
Then
\begin{align*}
\overline{d}\omega &  =-\overline{d}i_{\mathscr{E}}^{\ast}(d^{V}\theta)\\
&  =i_{\mathscr{E}}^{\ast}(d^{V}\overline{d}\theta)\\
&  =i_{\mathscr{E}}^{\ast}(d^{V}(d^{V}\mathscr{L}-\boldsymbol{E}%
(\mathscr{L})))\\
&  =-d^{V}i_{\mathscr{E}}^{\ast}(\boldsymbol{E}(\mathscr{L}))\\
&  =0.
\end{align*}
Since $\omega$ is $\overline{d}$-closed we may take its cohomology class
$\boldsymbol{\omega}:=[\omega]\in\boldsymbol{\Lambda}^{2}(\boldsymbol{P}%
)^{n-1}$. Now, $\boldsymbol{\omega}$ is canonical, as proved in what follows.

\begin{enumerate}
\item $\boldsymbol{\omega}$ does not depend on the choice of $\theta$. Indeed,
let $\theta^{\prime}:=\theta+\overline{d}\lambda$ be another Legendre form,
$\lambda\in\mathscr{C}\Lambda^{1}(J^{\infty})\otimes\overline{\Lambda}{}%
^{n-2}(J^{\infty})$ and $\omega^{\prime}:=-i_{\mathscr{E}}^{\ast}(d^{V}%
\theta^{\prime})$. Then $\omega^{\prime}=-i_{\mathscr{E}}^{\ast}(d^{V}%
\theta+d^{V}\overline{d}\lambda)=\omega+\overline{d}i_{\mathscr{E}}^{\ast
}(d^{V}\lambda)$, so that $[\omega]=[\omega^{\prime}]$.

\item $\boldsymbol{\omega}$ does not depend on the choice of $\mathscr{L}$.
Indeed, let $\mathscr{L}$ be a trivial Lagrangian density, i.e.,
$\mathscr{L}=\overline{d}\nu$ for some $\nu\in\overline{\Lambda}{}%
^{n-1}(J^{\infty})$. Then $\boldsymbol{S}=0$, $\boldsymbol{E}(\mathscr{L})=0$
and $d^{V}\mathscr{L}-\boldsymbol{E}(\mathscr{L})=-\overline{d}d^{V}\nu$. This
proves that $-d^{V}\nu$ is a Legendre form, so that $\boldsymbol{\omega
}=[i_{\mathscr{E}}^{\ast}(d^{V}d^{V}\nu)]=0$.
\end{enumerate}

Finally, $\boldsymbol{d\omega}=[d^{V}\omega]=0$.
\end{proof}

Notice that the above theorem can be generalized to the case of a Lagrangian field theory subject to constraints in the form of (the infinite prolongation of) a PDE $\mathscr{F}\subset J^\infty$, under suitable cohomological conditions on $\mathscr{F}$. Constrained Lagrangian theories will be considered somewhere else.

A general coordinate formula for $\omega$ may be found, for instance, in
\cite{lw90}. The expression of $\omega$ for specific Lagrangian theories may
be found, for instance, in \cite{c88,cw87,js02,lw90,r04}. However, we stress
that, in general, there is no distinguished representative $\omega$ in
$\boldsymbol{\omega}$.

\subsection{\textquotedblleft Symplectic Version\textquotedblright\ of I
Noether Theorem\label{Sec8}}

Let $(\pi,\boldsymbol{S})$ be a Lagrangian field theory and $\chi\in\varkappa$
the generating section of an higher symmetry of $\pi$. In view of isomorphism
$\varkappa\simeq\mathbf{D}(\boldsymbol{M})^{\bullet}$, $\chi$ may be
understood as a secondary vector field on $\boldsymbol{M}$. By definition,
$\chi$ is a \emph{Noether symmetry} of $(\pi,\boldsymbol{S})$ iff
$\mathcal{L}_{\chi}\boldsymbol{S}=0$, or, which is the same, $i_{\chi
}\boldsymbol{dS}=0$. In terms of a Lagrangian density $\mathscr{L}$ the last
equality reads as $i_{\rE_{\chi}}d^{V}\mathscr{L}=\overline{d}\sigma$ for some
$\sigma\in\overline{\Lambda}{}^{n-1}(J^{\infty})$. Using (\ref{Legendre}) one
gets $i_{\rE_{\chi}}(\boldsymbol{E}(\mathscr{L})+\overline{d}\theta
)=\overline{d}\sigma$. In view of isomorphism $\eta_{0}^{\ast}$, this implies
$\overline{d}(\sigma-i_{\rE_{\chi}}\theta)=\langle\boldsymbol{E}%
(\mathscr{L}),\chi\rangle$ and, pulling-back to $\mathscr{E}$,
\[
\overline{d}i_{\mathscr{E}}^{\ast}(\sigma-i_{\rE_{\chi}}\theta)=0.
\]
We have thus shown that $j:=i_{\mathscr{E}}^{\ast}(\sigma-i_{\rE_{\chi}}%
\theta)\in\overline{\Lambda}{}^{n-1}(\mathscr{E})$ is a conserved current of
$\mathscr{E}$ and this is, basically, the content of the \emph{first Noether
theorem}. Any such conserved current is called a \emph{Noether current} of
$(\pi,\boldsymbol{S})$. The associated conservation law $\boldsymbol{f}%
:=[j]\in\overline{H}{}^{n-1}(\mathscr{E})\subset\boldsymbol{C}^{\infty
}(\boldsymbol{P})^{\bullet}$ is called a \emph{Noether charge}. Notice that
neither $j$ nor $\boldsymbol{f}$ are uniquely determined by $\chi$ in general.

It is well known that if $\chi\in\varkappa$ is a Noether symmetry of
$(\pi,\boldsymbol{S})$, then $\chi|_{\mathscr{E}}\in\varkappa|_{\mathscr{E}}$
is the generating section of a symmetry of $\mathscr{E}$, i.e., $\ell
_{\boldsymbol{E}(\mathscr{L})}\chi|_{\mathscr{E}}=0$. This can be easily
proved by means of the following useful

\begin{lemma}
\label{Lemma}Let $\varphi\in\varkappa^{\dag}$,
$\mathscr{F}:=\mathscr{E}_{\varphi}^{(\infty)}\subset J^{\infty}$. For any
$\chi\in\varkappa$,
\[
(\mathcal{L}_{\rE_{\chi}}\varphi)|_{\mathscr{F}}=\ell_{\varphi}\overline{\chi
},\quad\overline{\chi}:=\chi|_{\mathscr{F}}%
\]
In particular, $(\mathcal{L}_{\rE_{\chi}}\varphi)|_{\mathscr{F}}\in
\varkappa|_{\mathscr{E}}^{\dag}$ and it does only depend on the values of
$\chi$ on $\mathscr{F}$.
\end{lemma}

\begin{proof}
For any $\chi_{1}\in\varkappa$, put $\overline{\chi}{}_{1}:=\chi
_{1}|_{\mathscr{F}}$. Similarly, for a (local) function $f$ on $J^{\infty}$,
put $\overline{f}:=f|_{\mathscr{F}}$. Compute%
\begin{align*}
i_{\rE_{\chi_{1}}}\mathcal{L}_{\rE_{\chi}}\varphi &  =i_{[\rE_{\chi_{1}%
},\rE_{\chi}]}\varphi+\mathcal{L}_{\rE_{\chi}}i_{\rE_{\chi_{1}}}\varphi\\
&  =i_{\rE_{\{\chi,\chi_{1}\}}}\varphi+\mathcal{L}_{\rE_{\chi}}\langle
\varphi,\chi_{1}\rangle\\
&  =\langle\varphi,\{\chi,\chi_{1}\}\rangle+\mathcal{L}_{\rE_{\chi}}%
\langle\varphi,\chi_{1}\rangle.
\end{align*}
Since $\varphi|_{\mathscr{F}}=0$, we have
\begin{equation}
(i_{\rE_{\chi_{1}}}\mathcal{L}_{\rE_{\chi}}\varphi)|_{\mathscr{F}}%
=(\mathcal{L}_{\rE_{\chi}}\langle\varphi,\chi_{1}\rangle)|_{\mathscr{F}}.
\label{Eq0}%
\end{equation}

Now, let $\varphi$, $\chi$ and $\chi_{1}$ be locally given by $\varphi
=\varphi_{\alpha}\partial^{\dag\alpha}$, $\chi=\chi^{\beta}\partial_{\beta}$,
$\chi_{1}=\chi_{1}^{\gamma}\partial_{\gamma}$, $\ldots,\varphi_{\alpha}%
,\ldots,\chi^{\beta},\ldots,\chi_{1}^{\gamma},\ldots$ local functions on
$J^{\infty}$. Then locally,
\[
\mathcal{L}_{\rE_{\chi}}\langle\varphi,\chi_{1}\rangle=D_{I}\chi^{\beta
}\partial_{\beta}^{I}(\varphi_{\alpha}\chi_{1}^{\alpha})\overline{\mathrm{d}%
}{}^{n}x=[D_{I}\chi^{\beta}(\partial_{\beta}^{I}\varphi_{\alpha})\chi
_{1}^{\alpha}+D_{I}\chi^{\beta}(\partial_{\beta}^{I}\chi_{1}^{\alpha}%
)\varphi_{\alpha}]\overline{\mathrm{d}}{}^{n}x.
\]
Since $\varphi_{\alpha}|_{\mathscr{F}}=0$, $\alpha=1,\ldots,m$, we have
locally%
\begin{equation}
(\mathcal{L}_{\rE_{\chi}}\langle\varphi,\chi_{1}\rangle)|_{\mathscr{F}}%
=D_{I}^{\mathscr{F}}\overline{\chi}{}^{\beta}\overline{(\partial_{\beta}%
^{I}\varphi_{\alpha})}\overline{\chi}{}{}_{1}^{\alpha}\overline{\mathrm{d}}%
{}^{n}x=\langle\ell_{\varphi}\overline{\chi},\overline{\chi}{}_{1}\rangle.
\label{Eq1}%
\end{equation}
Using (\ref{Eq1}) into (\ref{Eq0}) we get
\[
i_{\rE_{\chi_{1}}}^{\mathscr{F}}(\mathcal{L}_{\rE_{\chi}}\varphi
)|_{\mathscr{F}}=(i_{\rE_{\chi_{1}}}\mathcal{L}_{\rE_{\chi}}\varphi
)|_{\mathscr{F}}=\langle\ell_{\varphi}\overline{\chi},\overline{\chi}{}%
_{1}\rangle=i_{\rE_{\chi_{1}}}^{\mathscr{F}}\ell_{\varphi}\overline{\chi},
\]
where $i_{\rE_{\chi_{1}}}^{\mathscr{F}}$ is the restriction to $\mathscr{F}$
of the operator $i_{\rE_{\chi_{1}}}$ (see Section \ref{HorCalc}). From the
arbitrariness of $\chi_{1}$ the result follows.
\end{proof}

Now, let $\chi\in\varkappa$ be a Noether symmetry of the Lagrangian theory
$(\pi,\boldsymbol{S})$, and $\mathscr{L}$ a Lagrangian density. Then, in view
of Lemma \ref{Lemma},
\begin{equation}
\ell_{\boldsymbol{E}(\mathscr{L})}\chi|_{\mathscr{E}}=(\mathcal{L}_{\rE_{\chi
}}\boldsymbol{E}(\mathscr{L}))|_{\mathscr{E}}, \label{Eq4}%
\end{equation}
and
\[
\mathcal{L}_{\rE_{\chi}}\boldsymbol{E}(\mathscr{L})=\mathcal{L}_{\rE_{\chi}%
}(d^{V}\mathscr{L}-\overline{d}\theta)=d^{V}(i_{\rE_{\chi}}d^{V}%
\mathscr{L})+\overline{d}(\mathcal{L}_{\rE_{\chi}}\theta)=\overline
{d}(\mathcal{L}_{\rE_{\chi}}\theta-d^{V}\sigma).
\]
This shows that the horizontal cohomology class $[\mathcal{L}_{\rE_{\chi}%
}\boldsymbol{E}(\mathscr{L})]\in\boldsymbol{\Lambda}^{1}(\boldsymbol{M}%
)^{n}\simeq\varkappa^{\dag}$ is zero (and so is its \textquotedblleft
restriction\textquotedblright\ to $\mathscr{E}$) and, therefore,
$\ell_{\boldsymbol{E}(\mathscr{L})}\chi|_{\mathscr{E}}=0$ (see the final
comment in Example \ref{Ex1}).

The above remark proves that if $\chi$ is a Noether symmetry, then
$\boldsymbol{X}:=\chi|_{\mathscr{E}}\in\ker\ell_{\boldsymbol{E}(\mathscr{L})}%
\simeq\mathbf{D}(\boldsymbol{P})^{0}$ is a secondary vector field on
$\boldsymbol{P}$. Let $\boldsymbol{f}\in\mathbf{C}^{\infty}(\boldsymbol{P}%
)^{n-1}$ be, as above, a Noether charge associated to $\chi$.

\begin{proposition}
\label{I Noether}$\boldsymbol{df}=-i_{\boldsymbol{X}}\boldsymbol{\omega}$ (see
Equation 22 in \cite{lw90}).
\end{proposition}

\begin{proof}
Let $j$, $\sigma$, $\theta$ and $\omega$ be as above. First of all, notice that 
\begin{align*}
\overline{d}{}^\mathscr{E}(d^{V}\sigma-\mathcal{L}_{\rE_{\chi}}\theta)|_\mathscr{E} & =  \overline{d}(d^{V}\sigma-\mathcal{L}_{\rE_{\chi}}\theta)|_\mathscr{E} \\
                                                                                    & = -\mathcal{L}_{\rE_{\chi}}\boldsymbol{E}(\mathscr{L})|_\mathscr{E} \\
                                                                                    & = -\ell_{\boldsymbol{E}(\mathscr{L})}\chi|_{\mathscr{E}} \\
                                                                                    & =0
\end{align*}
and, therefore, $(d^{V}\sigma-\mathcal{L}_{\rE_{\chi}}\theta)|_\mathscr{E} = \overline{d}{}^\mathscr{E} \nu |_\mathscr{E} $ for some $\nu \in \mathscr{C}\Lambda^{1}(J^{\infty})\otimes\overline{\Lambda}{}^{n-2}(J^{\infty})$ and, in its turn, $i_{\mathscr{E}}^{\ast}(d^{V}\sigma-\mathcal{L}_{\rE_{\chi}}\theta)= \overline{d} i_{\mathscr{E}}^\ast (\nu)$. Then
\begin{align*}
\boldsymbol{df}  &  =[d^{V}j]\\
&  =[d^{V}i_{\mathscr{E}}^{\ast}(\sigma-i_{\rE_{\chi}}\theta)]\\
&  =[i_{\mathscr{E}}^{\ast}(d^{V}\sigma-d^{V}i_{\rE_{\chi}}\theta)]\\
&  =[i_{\mathscr{E}}^{\ast}(d^{V}\sigma-\mathcal{L}_{\rE_{\chi}}%
\theta+i_{\rE_{\chi}}d^{V}\theta)]\\
&  =[i_{\mathscr{E}}^{\ast}(d^{V}\sigma-\mathcal{L}_{\rE_{\chi}}%
\theta)+i_{\rE_{\chi}|_{\mathscr{E}}}i_{\mathscr{E}}^{\ast}(d^{V}\theta)]\\
&  =[i_{\mathscr{E}}^{\ast}(d^{V}\sigma-\mathcal{L}_{\rE_{\chi}}%
\theta)-i_{\rE_{\chi}|_{\mathscr{E}}}\omega]\\
&  =[\overline{d} i_{\mathscr{E}}^\ast (\nu)
]-i_{\boldsymbol{X}}\boldsymbol{\omega} \\
& = -i_{\boldsymbol{X}}\boldsymbol{\omega}.
\end{align*}
\end{proof}

Notice that Proposition \ref{I Noether} resembles very closely the analogous
result in Hamiltonian mechanics. Moreover, if $\mathscr{E}_{\boldsymbol{E}%
(\mathscr{L})}$ is an irreducible equation then $\boldsymbol{d}:\boldsymbol{C}%
^{\infty}(\boldsymbol{P})^{n-1}\longrightarrow\boldsymbol{\Lambda}%
^{1}(\boldsymbol{M})^{n-1}$ is injective \cite{kv98,v84} modulo obstructions
in $H^{n-1}(E)\subset\overline{H}{}^{n-1}(\mathscr{E})$. Thus,
$\boldsymbol{df}$ determines the \textquotedblleft\emph{non-trivial
conservation law}\textquotedblright\ (see \cite{b...99}) $\boldsymbol{f}%
+H^{n-1}(E)\in\overline{H}{}^{n-1}(\mathscr{E})/H^{n-1}(E)$ and is interpreted
as the \emph{generating section }of it. Proposition \ref{I Noether} can be
then understood as a way to compute the generating section of the non-trivial
conservation law associated to a Noether symmetry.

\subsection{\textquotedblleft Symplectic Version\textquotedblright\ of
(Infinitesimal) II Noether Theorem\label{Sec9}}

First of all, recall that the operator $\ell_{\boldsymbol{E}(\mathscr{L})}%
:\varkappa|_{\mathscr{E}}\longrightarrow\varkappa|_{\mathscr{E}}^{\dag}$ is
self-adjoint, i.e., $\ell_{\boldsymbol{E}(\mathscr{L})}=\ell_{\boldsymbol{E}%
(\mathscr{L})}^{\dag}:\varkappa|_{\mathscr{E}}\longrightarrow\varkappa
|_{\mathscr{E}}^{\dag}$. This fact is key in the calculus of variations
\cite{v84} and will be crucial in what follows (for a proof see, for instance,
\cite{b...99,v84} - see also \cite{a92} for an alternative approach).

The usual definition of (infinitesimal) gauge symmetries of a Lagrangian field
theory is the following (see \cite{lw90}).

\begin{definition}
A $\emph{Noether}$ \emph{gauge} (or \emph{local}) \emph{symmetry} of the
Lagrangian theory $(\pi,\boldsymbol{S})$ is a horizontal linear differential
operator $G:Q\longrightarrow\varkappa$ such that $G(\varepsilon)$ is a Noether
symmetry for any $\varepsilon\in Q$.
\end{definition}

We added the prefix \textquotedblleft Noether\textquotedblright\ in the above
definition of a \textquotedblleft gauge symmetry\textquotedblright\ to
distinguish it from an alternative (and, generally, inequivalent) definition
that will be proposed below. Physicists say sometimes that $G$ is a
\emph{Noether symmetry depending on the arbitrary parameters }$\varepsilon$.

The second Noether theorem states that, in presence of gauge symmetries, there
are relations among the Euler-Lagrange equations. Namely, for all
$\varepsilon\in Q$,
\begin{align*}
0  &  =\mathcal{L}_{G(\varepsilon)}\boldsymbol{S}\\
&  =i_{G(\varepsilon)}\boldsymbol{dS}\\
&  =\int i_{\rE_{G(\varepsilon)}}d^{V}\mathscr{L}\\
&  =\int\tilde{\ell}{}_{\mathscr{L}}(G(\varepsilon))\\
&  =\int\langle1,(\tilde{\ell}{}_{\mathscr{L}}\circ G)(\varepsilon)\rangle\\
&  =\int\langle(\tilde{\ell}{}_{\mathscr{L}}\circ G)^{\dag}(1),\varepsilon
\rangle\\
&  =\int\langle G^{\dag}(\tilde{\ell}{}_{\mathscr{L}}^{\dag}1),\varepsilon
\rangle\\
&  =\int\langle G^{\dag}(\boldsymbol{E}(\mathscr{L})),\varepsilon\rangle,
\end{align*}
and it follows from the arbitrariness of $\varepsilon$ that $G^{\dag
}(\boldsymbol{E}(\mathscr{L}))=0$. These relations are traditionally called
\emph{Noether identities}.

An \textquotedblleft infinitesimal version\textquotedblright\ of the second
Noether theorem can be formulated. First of all, notice that, since
$G(\varepsilon)$ is a Noether symmetry (so that $G(\varepsilon)|_{\mathscr{E}}%
$ is the generating section of a symmetry of $\mathscr{E}$) for all
$\varepsilon$, one also has $0=\ell_{\boldsymbol{E}(\mathscr{L})}%
G(\varepsilon)|_{\mathscr{E}}=(\ell_{\boldsymbol{E}(\mathscr{L})}\circ
G^{\mathscr{E}})(\varepsilon|_{\mathscr{E}})$ and, from the arbitrariness of
$\varepsilon$,
\begin{equation}
\ell_{\boldsymbol{E}(\mathscr{L})}\circ G^{\mathscr{E}}=0. \label{Eq2}%
\end{equation}
Identity (\ref{Eq2}) may be interpreted by saying that the linearized
Euler-Lagrange equations admit \textquotedblleft\emph{gauge symmetries}%
\textquotedblright. Indeed, if $\chi\in\varkappa|_{\mathscr{E}}$ is in the
kernel of $\ell_{\boldsymbol{E}(\mathscr{L})}$ so is the \textquotedblleft%
\emph{gauge transformed}\textquotedblright\ element $\chi+G^{\mathscr{E}}%
(\epsilon)$, for any arbitrary $\epsilon\in Q|_{\mathscr{E}}$. In particular,
the linearized Euler-Lagrange equations are, in a sense, \textquotedblleft
underdetermined\textquotedblright.

By passing to the adjoint operators in (\ref{Eq2}) and using the
self-adjointness of $\ell_{\boldsymbol{E}(\mathscr{L})}$ we get
\begin{equation}
(G^{\mathscr{E}})^{\dag}\circ\ell_{\boldsymbol{E}(\mathscr{L})}=0. \label{Eq3}%
\end{equation}
This shows that there are relations among the linearized Euler-Lagrange
equations and that they are, in a sense, \textquotedblleft
constrained\textquotedblright. Thus, \textquotedblleft\emph{infinitesimal
gauge symmetries correspond to infinitesimal constraints}\textquotedblright%
\ via adjunction \cite{lw90}. Identities (\ref{Eq3}) (and sometimes the
operator $(G^{\mathscr{E}})^{\dag}$ itself) are called \emph{infinitesimal
Noether identities}.

Now let $\Delta_{1}:\varkappa|_{\mathscr{E}}^{\dag}\longrightarrow P_{2}$ be a
compatibility operator for $\ell_{\boldsymbol{E}(\mathscr{L})}$. Consider also
the adjoint operator $\Delta_{1}^{\dag}:P_{2}^{\dag}\longrightarrow
\varkappa|_{\mathscr{E}}$. In particular, $\Delta_{1}\circ
\ell_{\boldsymbol{E}(\mathscr{L})}=0$ and (using again the self-adjointness
of $\ell_{\boldsymbol{E}(\mathscr{L})}$) $\ell_{\boldsymbol{E}(\mathscr{L})}%
\circ\Delta_{1}^{\dag}=0$. In view of the last identity, if $\chi\in
\varkappa|_{\mathscr{E}}$ is in the kernel of $\ell_{\boldsymbol{E}%
(\mathscr{L})}$ so is the\ element $\chi+\Delta_{1}^{\dag}\vartheta$, for any
arbitrary $\vartheta\in P_{2}^{\dag}$. Notice also that, in view of
Proposition \ref{Prop1}, all infinitesimal Noether identities
$(G^{\mathscr{E}})^{\dag}$ \textquotedblleft are generated by $\Delta_{1}%
$\textquotedblright\ in the sense that $(G^{\mathscr{E}})^{\dag}=\nabla
\circ\Delta_{1}$ for some horizontal differential operator $\nabla
:P_{2}\longrightarrow Q|_{\mathscr{E}}^\dag$. Similarly, by passing to the adjoint
operators, we see that all \emph{infinitesimal gauge symmetries}
$G^{\mathscr{E}}$ are generated by $\Delta_{1}^{\dag}$, i.e., $G^{\mathscr{E}}%
=\Delta_{1}^{\dag}\circ\nabla^{\dag}$ for some $\nabla^{\dag}:Q|_{\mathscr{E}}%
\longrightarrow P_{2}^{\dag}$. These simple remarks suggest a more natural
definition of infinitesimal gauge symmetries.

\begin{definition}
\label{DefGau}A \emph{gauge symmetry} of the Lagrangian theory $(\pi
,\boldsymbol{S})$ is an element in the image of the adjoint operator $\Delta_{1}^{\dag
}$ of a compatibility operator $\Delta_{1}$ for $\ell_{\boldsymbol{E}%
(\mathscr{L})}$.
\end{definition}

We will sometimes denote by $\mathfrak{g}:=\operatorname{im}\Delta_{1}^{\dag}$
the set of gauge symmetries. Notice that, in view of Theorem \ref{SpencTh},
the above definition is independent of the choice of $\Delta_{1}$. Moreover,
while it is clear that $\operatorname{im}G^{\mathscr{E}}\subset\mathfrak{g}$
for any Noether gauge symmetry $G$, to the author knowledge it has not been
determined yet in full rigour and generality if $\mathfrak{g}$ is generated by
the images of Noether gauge symmetries or not. Therefore, we prefer to adopt
definition \ref{DefGau}. This choice is strengthened even more by the results
presented in the remaining part of this section.

Consider the natural $\mathbb{R}$-linear map
\[
\boldsymbol{\Omega}:\mathbf{D}(\boldsymbol{P})^{\bullet}\ni\boldsymbol{X}%
\longmapsto\boldsymbol{\Omega}(\boldsymbol{X}):=i_{\boldsymbol{X}%
}\boldsymbol{\omega}\in\boldsymbol{\Lambda}^{1}(\boldsymbol{P})^{\bullet}.
\]

\begin{definition}
The kernel $\ker\boldsymbol{\Omega}\subset\mathbf{D}(\boldsymbol{P}%
)^{\bullet}$ is called the \emph{degeneracy distribution} of
$\boldsymbol{\omega}$ and will be also denoted by $\ker\boldsymbol{\omega}$. The secondary $2$-form $\boldsymbol{\omega}$ is said to be
1) \emph{weakly symplectic} (or \emph{non-degenerate}) iff $\ker\boldsymbol{\omega}=0$, 2) \emph{strongly symplectic} (or, simply, \emph{symplectic}) iff $\boldsymbol{\Omega}$ is
an isomorphism. 
\end{definition}

In order to better characterize $\boldsymbol{\omega}$ it is desirable to
describe its degeneracy distribution.

First of all, notice that, since $\boldsymbol{\omega}$ is closed,
$\ker\boldsymbol{\omega}$ \emph{is a secondary involutive distribution}, i.e.,
it is a graded Lie subalgebra in $\mathbf{D}(\boldsymbol{P})^{\bullet}$.
Indeed, let $\boldsymbol{X},\boldsymbol{Y}\in\ker\boldsymbol{\omega}$ then
\[
\boldsymbol{\Omega}([\boldsymbol{X},\boldsymbol{Y}])=i_{[\boldsymbol{X}%
,\boldsymbol{Y}]}\boldsymbol{\omega}=[i_{\boldsymbol{X}},\mathcal{L}%
_{\boldsymbol{Y}}]\boldsymbol{\omega}=[i_{\boldsymbol{X}},[i_{\boldsymbol{Y}%
},\boldsymbol{d}]]\boldsymbol{\omega}=0,
\]
i.e., $[\boldsymbol{X},\boldsymbol{Y}]\in\ker\boldsymbol{\omega}$.

Denote by $\boldsymbol{\Omega}^{r}:\mathbf{D}(\boldsymbol{P})^{r}%
\longrightarrow\boldsymbol{\Lambda}^{1}(\boldsymbol{P})^{r+n-1}$ the
restriction of $\boldsymbol{\Omega}$ to $\mathbf{D}(\boldsymbol{P})^{r}$,
$r=0,\ldots,n$. Obviously, $\boldsymbol{\Omega}^{r}=0$ for $r>1$,
independently of the Lagrangian theory. For this reason, every degree $r>1$
secondary vector field over $\boldsymbol{P}$ is said to be a \emph{trivial
element in }$\ker\boldsymbol{\omega}$. Thus, non-trivial elements in
$\ker\boldsymbol{\omega}$ must be searched in $\mathbf{D}(\boldsymbol{P})^{0}$
and $\mathbf{D}(\boldsymbol{P})^{1}$. In the following we will
\textquotedblleft describe\textquotedblright\ such elements. Put
$\operatorname{tker}\boldsymbol{\omega}:=\ker\boldsymbol{\omega\cap}%
\bigoplus_{r>1}\mathbf{D}(\boldsymbol{P})^{r}$.

\begin{theorem}
\label{KerOm}Diagrams
\begin{equation}%
\begin{array}
[c]{c}%
\xymatrix{  & & \mathbf{D}(\boldsymbol{P})^0 \ar[r]^-{\boldsymbol{\Omega}^0} &    \boldsymbol{\Lambda}^1(\boldsymbol{P})^{n-1} \ar@{}@<0.2ex>[d]^{\begin{sideways}$\widetilde{\quad\quad}$\end{sideways}} \ar[d]                 &          \\
0 \ar[r]& \operatorname{im}\Delta_1^\dag\, \ar@{^{(}->}[r] &  \ker \ell_{\boldsymbol{E}(\mathscr{L})} \ar@{->>}[r] \ar@{}@<-1.2ex>[u]^{\begin{sideways}$\widetilde{\quad\quad}$\end{sideways}} \ar[u]& \ker \ell_{\boldsymbol{E}(\mathscr{L})}/\operatorname{im}\Delta_1^\dag \ar[r] & 0}
\end{array}
\label{KerOm1}%
\end{equation}
and%
\begin{equation}%
\begin{array}
[c]{c}%
\xymatrix{  &  \mathbf{D}(\boldsymbol{P})^1 \ar[r]^-{\boldsymbol{\Omega}^1} &    \boldsymbol{\Lambda}^1(\boldsymbol{P})^{n} \ar@{}@<0.2ex>[d]^{\begin{sideways}$\widetilde{\quad\quad}$\end{sideways}} \ar[d]    & &                       \\
0 \ar[r]& \ker\Delta_1 / \operatorname{im} \ell_{\boldsymbol{E}(\mathscr{L})}\, \ar@{^{(}->}[r] \ar[u] \ar@{}@<-1.2ex>[u]^{\begin{sideways}$\widetilde{\quad\quad}$\end{sideways}}&  \operatorname{coker} \ell_{\boldsymbol{E}(\mathscr{L})} \ar@{->>}[r]  & \varkappa |^\dag_{\mathscr{E}}/\operatorname{ker}\Delta_1 \ar[r] & 0}
\end{array}
\label{KerOm2}%
\end{equation}
commute.
\end{theorem}

\begin{proof}
The vertical arrows in Diagram (\ref{KerOm1}) are described in Section
\ref{SpencGold}. Thus, let $\boldsymbol{X}=\rE_{\chi}\in$ $\mathbf{D}%
(\boldsymbol{P})^{0}$, $\chi\in\varkappa|_{\mathscr{E}}$, $\ell
_{\boldsymbol{E}(\mathscr{L})}\chi=0$. Let $\tilde{\chi}\in\varkappa$ be such
that $\tilde{\chi}|_{\mathscr{E}}=\chi$. Now, $\boldsymbol{\Omega}%
^{0}(\boldsymbol{X})=i_{\boldsymbol{X}}\boldsymbol{\omega}=[i_{\mathscr{E}}%
^{\ast}(-i_{\rE_{\tilde{\chi}}}d^{V}\theta)]$, $\theta$ being a Legendre form.
Put $\tilde{\square}:=(\eta_{0}^{\ast})^{-1}(-i_{\rE_{\tilde{\chi}}}%
d^{V}\theta)\in\mathscr{C}\mathrm{Diff}(\varkappa,\overline{\Lambda}{}%
^{n-1}(J^{\infty}))$ and $\square:=\tilde{\square}{}^{\mathscr{E}}%
\in\mathscr{C}\mathrm{Diff}(\varkappa|_{\mathscr{E}},\overline{\Lambda}%
{}^{n-1}(\mathscr{E}))$. Then, obviously, $\eta_{\boldsymbol{E}(\mathscr{L})}%
^{\ast}(\square)=i_{\mathscr{E}}^{\ast}(-i_{\rE_{\tilde{\chi}}}d^{V}\theta)$.
Show that $\overline{d}\circ\square=\Delta_{\chi}\circ\ell_{\boldsymbol{E}%
(\mathscr{L})}$ where $\Delta_{\chi}\in\mathscr{C}\mathrm{Diff}(\varkappa
|_{\mathscr{E}}^{\dag},\overline{\Lambda}{}^{n}(\mathscr{E}))$ is defined by
putting $\Delta_{\chi}\varphi:=\langle\varphi,\chi\rangle$, $\varphi
\in\varkappa|_{\mathscr{E}}^{\dag}$ (thus, $\Delta_{\chi}$ is actually a
$C^{\infty}(\mathscr{E})$-linear map). Indeed, let $\chi_{1}\in\varkappa$ and
put $\overline{\chi}_{1}:=\chi_{1}|_{\mathscr{E}}$. Compute
\begin{align*}
(\overline{d}\circ\square)(\overline{\chi}_{1})  &  =\overline{d}%
(\tilde{\square}\chi_{1})|_{\mathscr{E}}\\
&  =\overline{d}(-i_{\rE_{\chi_{1}}}i_{\rE_{\tilde{\chi}}}d^{V}\theta
)|_{\mathscr{E}}\\
&  =(-i_{\rE_{\tilde{\chi}}}i_{\rE_{\chi_{1}}}d^{V}\overline{d}\theta
)|_{\mathscr{E}}\\
&  =(i_{\rE_{\tilde{\chi}}}i_{\rE_{\chi_{1}}}d^{V}\boldsymbol{E}%
(\mathscr{L}))|_{\mathscr{E}}\\
&  =(i_{\rE_{\tilde{\chi}}}\mathcal{L}_{\rE_{\chi_{1}}}\boldsymbol{E}%
(\mathscr{L}))|_{\mathscr{E}}-(i_{\rE_{\tilde{\chi}}}d^{V}\langle
\boldsymbol{E}(\mathscr{L}),\chi_{1}\rangle)|_{\mathscr{E}}\\
&  =\langle\ell_{\boldsymbol{E}(\mathscr{L})}\overline{\chi}_{1},\chi
\rangle-(\mathcal{L}_{\rE_{\tilde{\chi}}}\langle\boldsymbol{E}%
(\mathscr{L}),\chi_{1}\rangle)|_{\mathscr{E}}\\
&  =\langle\ell_{\boldsymbol{E}(\mathscr{L})}\overline{\chi}_{1},\chi
\rangle-\langle\ell_{\boldsymbol{E}(\mathscr{L})}\chi,\overline{\chi}%
_{1}\rangle\\
&  =(\Delta_{\chi}\circ\ell_{\boldsymbol{E}(\mathscr{L})})(\overline{\chi}%
_{1}),
\end{align*}
where we used Identities (\ref{Eq1}) and (\ref{Eq4}). It follows from the
arbitrariness of $\chi_{1}$ that $\overline{d}\circ\square=\Delta_{\chi}%
\circ\ell_{\boldsymbol{E}(\mathscr{L})}$. Therefore, $i_{\boldsymbol{X}%
}\boldsymbol{\omega}$ corresponds to $\Delta_{\chi}^{\dag}1+\operatorname{im}%
\Delta_{1}^{\dag}\in\ker\ell_{\boldsymbol{E}(\mathscr{L})}/\operatorname{im}%
\Delta_{1}^{\dag}$ via isomorphism $\boldsymbol{\Lambda}^{1}(\boldsymbol{P}%
)^{n-1}\simeq\ker\ell_{\boldsymbol{E}(\mathscr{L})}/\operatorname{im}%
\Delta_{1}^{\dag}$. Finally, it is easy to see that $\Delta_{\chi}^{\dag
}1=\chi$.

Now, consider diagram (\ref{KerOm2}) whose vertical arrows are described in
Section \ref{SpencGold} as well. Let $\varphi\in\varkappa|_{\mathscr{E}}%
^{\dag}$ and $j\in\overline{J}{}^{\infty}\varkappa$ be such that $\Delta
_{1}\varphi=0$ and $\overline{j}{}_{\infty}\varphi=h_{\boldsymbol{E}%
(\mathscr{L})}(j|_{\mathscr{E}})\in\overline{J}{}^{\infty}\varkappa
|_{\mathscr{E}}^{\dag}$. Since $\overline{J}{}^{\infty}\varkappa$ is
pro-finitely generated by elements of the form $\overline{j}%
{}_{\infty}\chi$,\footnote{This means that an element in $\overline{J}%
{}^{\infty}\varkappa$ may be understood as a(n equivalence class of) formal
infinite linear combination(s) of elements of the form $\overline{j}_{\infty}%
\chi$, $\chi\in\varkappa$. Notice that, in any case, all the following
computations remain still valid.}\label{page} $\chi\in\varkappa$, then $j=\sum
f\overline{j}{}_{\infty}\chi$ for some (generally, infinite in number)
$\ldots,f,\ldots\in C^{\infty}(J^{\infty})$ and $\ldots,\chi,\ldots
\in\varkappa$. Consequently, $\varphi=\sum f|_{\mathscr{E}}\ell
_{\boldsymbol{E}(\mathscr{L})}\chi|_{\mathscr{E}}$. Put $Z:=(\overline{S}%
\circ\eta_{0}^{-1})(j)=\sum\rE_{\chi}\otimes\overline{d}f\in V\mathrm{D}%
(J^{\infty})\otimes\overline{\Lambda}{}^{1}(J^{\infty})$ and recall that 1)
$Z$ restricts to $\mathscr{E}$ and 2) $\varphi+\operatorname{im}%
\ell_{\boldsymbol{E}(\mathscr{L})}\in\ker\Delta_{1}/\operatorname{im}%
\ell_{\boldsymbol{E}(\mathscr{L})}$ corresponds to $\boldsymbol{Z}%
:=[\overline{Z}]\in\mathbf{D}(\boldsymbol{P})^{1}$, $\overline{Z}\in
V\mathrm{D}(\mathscr{E})\otimes\overline{\Lambda}{}^{1}(\mathscr{E})$ being
the restriction of $Z$ to $\mathscr{E}$, via isomorphism $\ker\Delta
_{1}/\operatorname{im}\ell_{\boldsymbol{E}(\mathscr{L})}\simeq\mathbf{D}%
(\boldsymbol{P})^{1}$. Now, $\boldsymbol{\Omega}^{1}(\boldsymbol{Z}%
)=i_{\boldsymbol{Z}}\boldsymbol{\omega}=[i_{\overline{Z}}i_{\mathscr{E}}%
^{\ast}(d^{V}\theta)]=[i_{\mathscr{E}}^{\ast}(i_{Z}d^{V}\theta)]$. Compute
\begin{align*}
i_{Z}d^{V}\theta &  =\sum\overline{d}f\wedge i_{\rE_{\chi}}d^{V}\theta\\
&  =\overline{d}\rho-\sum f\overline{d}i_{\rE_{\chi}}d^{V}\theta\\
&  =\overline{d}\rho-\sum fi_{\rE_{\chi}}d^{V}\overline{d}\theta\\
&  =\overline{d}\rho+\sum(f\mathcal{L}_{\rE_{\chi}}\boldsymbol{E}%
(\mathscr{L})-fd^{V}\langle\boldsymbol{E}(\mathscr{L}),\chi\rangle),
\end{align*}
where $\rho=\sum fi_{\rE_{\chi}}d^{V}\theta\in\mathscr{C}\Lambda^{1}%
(J^{\infty})\otimes\overline{\Lambda}{}^{n-1}(J^{\infty})$. Therefore,
\begin{align*}
i_{\overline{Z}}i_{\mathscr{E}}^{\ast}(d^{V}\theta)  &  =\overline
{d}i_{\mathscr{E}}^{\ast}(\rho)+\sum i_{\mathscr{E}}^{\ast}(f\mathcal{L}%
_{\rE_{\chi}}\boldsymbol{E}(\mathscr{L})+fd^{V}\langle\boldsymbol{E}%
(\mathscr{L}),\chi\rangle)\\
&  =\overline{d}i_{\mathscr{E}}^{\ast}(\rho)+\sum\eta_{\boldsymbol{E}%
(\mathscr{L})}^{\ast}(f|_{\mathscr{E}}\ell_{\boldsymbol{E}(\mathscr{L})}%
\chi|_{\mathscr{E}})+\sum fd^{V}i_{\mathscr{E}}^{\ast}\langle\boldsymbol{E}%
(\mathscr{L}),\chi\rangle\\
&  =\overline{d}i_{\mathscr{E}}^{\ast}(\rho)+\eta_{\boldsymbol{E}%
(\mathscr{L})}^{\ast}(\varphi).
\end{align*}
Finally, $\boldsymbol{\Omega}^{1}(\boldsymbol{Z})=[\eta_{\boldsymbol{E}%
(\mathscr{L})}^{\ast}(\varphi)]$ corresponds to $\varphi^{\dag}1+\operatorname{im} \ell_{\boldsymbol{E}(\mathscr{L})} 
\in\operatorname{coker}\ell_{\boldsymbol{E}(\mathscr{L})}$ via isomorphism $\boldsymbol{\Lambda}%
^{1}(\boldsymbol{P})^{n}\simeq\operatorname{coker}\ell_{\boldsymbol{E}%
(\mathscr{L})}$. It is easy to prove that $\varphi^{\dag}1=\varphi$ and this
concludes the proof.
\end{proof}

Some corollaries are in order.

\begin{corollary}
There is a natural isomorphism $\ker\boldsymbol{\omega}\simeq\mathfrak{g}%
\oplus\operatorname{tker}\boldsymbol{\omega}$.
\end{corollary}

\begin{corollary}
\label{CorsubLie}$\mathfrak{g}\subset\ker\ell_{\boldsymbol{E}(\mathscr{L})}$
is a Lie subalgebra (see, for instance, \cite{bc08}).
\end{corollary}

\begin{corollary}
\label{CorNoethGau}Let $G:Q\longrightarrow\varkappa$ be a Noether gauge
symmetry. Then $\operatorname{im}G^{\mathscr{E}}\subset\ker\boldsymbol{\omega
}$ (see also \cite{lw90}).
\end{corollary}

\begin{corollary}
The secondary 2-form $\boldsymbol{\omega}$ is weakly symplectic iff it is strongly symplectic iff the
Euler-Lagrange equations $\mathscr{E}_{\boldsymbol{E}(\mathscr{L})}$ are irreducible.
\end{corollary}

\begin{proof}
In view of Theorem \ref{KerOm}, $\boldsymbol{\Omega}^{0}$ and
$\boldsymbol{\Omega}^{1}$ are isomorphisms iff $\mathscr{E}_{\boldsymbol{E}%
(\mathscr{L})}$ is an irreducible PDE (see Example \ref{IrrEq}). In view of the
$2$-lines Theorem \ref{slines} ($s=2$), irreducibility of
$\mathscr{E}_{\boldsymbol{E}(\mathscr{L})}$ implies, in its turn, that
$\operatorname{tker}\boldsymbol{\omega}=0$.
\end{proof}

\subsection{Gauge Invariant Secondary Functions\label{Sec10}}

Let $N$ be a smooth manifold and $\sigma\in\Lambda^{2}(N)$ a presymplectic
structure on it. There is no Poisson structure on $N$ associated to $\sigma$ .
However, a Poisson bracket may be introduced among \textquotedblleft gauge
invariant\textquotedblright\ functions on $N$, i.e., functions which are
constant along the leaves of the degeneracy distribution of $\sigma$. This is
precisely the Poisson bracket on the symplectic reduction of $(N,\sigma)$. In
this section we describe \textquotedblleft gauge invariant secondary
functions\textquotedblright\ on the CPS $\boldsymbol{P}$ and show that,
similarly to the standard situation, $\boldsymbol{\omega}$ induces a Lie bracket among them. Thus, the results presented in this
section are propaedeutic to a \textquotedblleft secondary symplectic
reduction\textquotedblright\ of $(\boldsymbol{P},\boldsymbol{\omega})$ (see
next section).

\begin{definition}
A secondary function $\boldsymbol{f}\in\boldsymbol{C}^{\infty}(\boldsymbol{P}%
)^{\bullet}$ is called \emph{gauge invariant} iff $\mathcal{L}_{\boldsymbol{Y}%
}\boldsymbol{f}=0$ for all $\boldsymbol{Y}\in\ker\boldsymbol{\omega}$.
\end{definition}

Let us describe gauge invariant elements in $\boldsymbol{C}^{\infty
}(\boldsymbol{P})^{n-1}$ and $\boldsymbol{C}^{\infty}(\boldsymbol{P})^{n}$.

\begin{proposition}
\label{PropGauInv}Any element in $\boldsymbol{C}^{\infty}(\boldsymbol{P}%
)^{n-1}$ is gauge invariant.
\end{proposition}

\begin{proof}
Recall that the map $\boldsymbol{\Omega}^{0}:\mathbf{D}(\boldsymbol{P}%
)^{0}\longrightarrow\mathbf{\Lambda}^{1}(\boldsymbol{P})^{n-1}$ is surjective
(see Theorem \ref{KerOm}). For any $\boldsymbol{f}\in\boldsymbol{C}^{\infty
}(\boldsymbol{P})^{n-1}$, let $\boldsymbol{X}\in\mathbf{D}(\boldsymbol{P}%
)^{0}$ be such that $\boldsymbol{\Omega}(\boldsymbol{X})=\boldsymbol{df}%
\in\mathbf{\Lambda}^{1}(\boldsymbol{P})^{n-1}$ and $\boldsymbol{Y}\in
\ker\boldsymbol{\omega}$. Then $\mathcal{L}_{\boldsymbol{Y}}\boldsymbol{f}%
=i_{\boldsymbol{Y}}\boldsymbol{df}=i_{\boldsymbol{Y}}i_{\boldsymbol{X}%
}\boldsymbol{\omega}=-i_{\boldsymbol{X}}i_{\boldsymbol{Y}}\boldsymbol{\omega
}=0$.
\end{proof}

Now, let $\boldsymbol{f}_{1},\boldsymbol{f}_{2}\in\boldsymbol{C}^{\infty
}(\boldsymbol{P})^{n-1}$ and $\boldsymbol{X}_{1},\boldsymbol{X}_{2}%
\in\mathbf{D}(\boldsymbol{P})^{0}$ be such that $\boldsymbol{\Omega
}(\boldsymbol{X}_{1})=\boldsymbol{df}_{1}$ and $\boldsymbol{\Omega
}(\boldsymbol{X}_{2})=\boldsymbol{df}_{2}$. Put $\{\boldsymbol{f}%
_{1},\boldsymbol{f}_{2}\}:=-i_{\boldsymbol{X}_{1}}i_{\boldsymbol{X}_{2}%
}\boldsymbol{\omega}\in\boldsymbol{C}^{\infty}(\boldsymbol{P})^{n-1}$.

\begin{corollary}
$(\boldsymbol{C}^{\infty}(\boldsymbol{P})^{n-1},\{{}\cdot{},{}\cdot{}\})$ is a
well defined Lie algebra.
\end{corollary}

\begin{proof}
In view of Proposition \ref{PropGauInv}, $\{\boldsymbol{f}_{1},\boldsymbol{f}%
_{2}\}$ is well defined for all $\boldsymbol{f}_{1},\boldsymbol{f}_{2}%
\in\boldsymbol{C}^{\infty}(\boldsymbol{P})^{n-1}$, i.e., it is independent of
the choice of $\boldsymbol{X}_{1},\boldsymbol{X}_{2}$. Skew-symmetry and the
Leibnitz rule follow (as in standard presymplectic geometry) from
$\boldsymbol{d\omega}=0$ and the fact that, if $\boldsymbol{\Omega
}(\boldsymbol{X}_{1})=\boldsymbol{df}_{1}$ and $\boldsymbol{\Omega
}(\boldsymbol{X}_{2})=\boldsymbol{df}_{2}$, then $\boldsymbol{\Omega
}([\boldsymbol{X}_{1},\boldsymbol{X}_{2}])=\boldsymbol{d}\{\boldsymbol{f}%
_{1},\boldsymbol{f}_{2}\}$.
\end{proof}

Notice that the existence of a natural Lie bracket among conservation laws of
an Euler-Lagrange equation was already known and may be also proved by
\textquotedblleft off shell\textquotedblright\ methods such as BRST ones (see,
for instance, \cite{bbh95,bbh00}).

\begin{proposition}
\label{GauInv}An element $\boldsymbol{F}\in\boldsymbol{C}^{\infty
}(\boldsymbol{P})^{n}$ is gauge invariant iff $\boldsymbol{dF}\in
\operatorname{im}\boldsymbol{\Omega}$.
\end{proposition}

\begin{proof}
If $\boldsymbol{dF}=\boldsymbol{\Omega}(\boldsymbol{Z})$ for some
$\boldsymbol{Z}\in\mathbf{D}(\boldsymbol{P})^{1}$, then $\boldsymbol{F}$ is
gauge invariant (see the proof of Proposition \ref{PropGauInv}). Vice versa,
suppose $\mathcal{L}_{\boldsymbol{Y}}\boldsymbol{F}=0$ for all $\boldsymbol{Y}%
\in\ker\boldsymbol{\omega}$. Let $\boldsymbol{F}=\int\rho$, $\rho\in
\overline{\Lambda}{}^{n}(\mathscr{E})$. Recall that $\boldsymbol{dF}%
=[d^{V}\rho]\in\boldsymbol{\Lambda}^{1}(\boldsymbol{P})^{n}$ corresponds to
$\square^{\dag}1+\operatorname{im}\ell_{\boldsymbol{E}(\mathscr{L})}%
\in\operatorname{coker}\ell_{\boldsymbol{E}(\mathscr{L})}$ via the isomorphism
$\boldsymbol{\Lambda}^{1}(\boldsymbol{P})^{n}\simeq\operatorname{coker}%
\ell_{\boldsymbol{E}(\mathscr{L})}$, $\square:\varkappa|_{\mathscr{E}}%
\longrightarrow\overline{\Lambda}{}^{n}(\mathscr{E})$ being any horizontal
differential operator such that $\eta_{\mathscr{E}}^{\ast}(\square)=d^{V}\rho
$. In view of Theorem \ref{KerOm}, $\boldsymbol{dF}\in\operatorname{im}%
\boldsymbol{\Omega}$ if $\Delta_{1}(\square^{\dag}1)=0$. Let $\chi=\Delta
_{1}^{\dag}\vartheta\in\varkappa|_{\mathscr{E}}$, $\vartheta\in P_{1}^{\dag}$,
and $\boldsymbol{Y}=\rE_{\chi}$. Then
\begin{align*}
0  &  =\mathcal{L}_{\boldsymbol{Y}}\boldsymbol{F}\\
&  =\int\mathcal{L}_{\rE_{\chi}}\rho\\
&  =\int i_{\rE_{\chi}}d^{V}\rho\\
&  =\int\square\chi\\
&  =\int(\square\circ\Delta_{1}^{\dag})(\vartheta)\\
&  =\int\langle(\square\circ\Delta_{1}^{\dag})(\vartheta),1\rangle\\
&  =\int\langle\vartheta,(\Delta_{1}\circ\square^{\dag})(1)\rangle.
\end{align*}
It follows from the arbitrariness of $\vartheta$ that $\Delta_{1}%
(\square^{\dag}1)=0$, and this concludes the proof.
\end{proof}

The Lie algebra $(\boldsymbol{C}^{\infty}(\boldsymbol{P})^{n-1},\{{}\cdot{}%
,{}\cdot{}\})$ acts naturally on gauge invariant elements in $\boldsymbol{C}%
^{\infty}(\boldsymbol{P})^{n}$. Indeed, let $\boldsymbol{F}\in\boldsymbol{C}%
^{\infty}(\boldsymbol{P})^{n}$ be a gauge invariant element and
$\boldsymbol{f}\in\boldsymbol{C}^{\infty}(\boldsymbol{P})^{n-1}$. Put
$\{\boldsymbol{f},\boldsymbol{F}\}:=\mathcal{L}_{\boldsymbol{X}}%
\boldsymbol{F}\in\boldsymbol{C}^{\infty}(\boldsymbol{P})^{n}$, $\boldsymbol{X}%
\in\mathbf{D}(\boldsymbol{P})^{0}$ being any secondary vector field such that
$\boldsymbol{\Omega}^{0}(\boldsymbol{X})=\boldsymbol{df}$. Exactly as above,
$\{\boldsymbol{f},\boldsymbol{F}\}$ is well defined. Moreover, it holds the

\begin{proposition}
$\{\boldsymbol{f},\boldsymbol{F}\}$ is gauge invariant.
\end{proposition}

\begin{proof}
Recall that, in view of Proposition \ref{GauInv}, $\boldsymbol{dF}%
\in\operatorname{im}\boldsymbol{\Omega}^{1}$, i.e., $\boldsymbol{dF}%
=i_{\boldsymbol{Z}}\boldsymbol{\omega}$ for some $\boldsymbol{Z}\in
\mathbf{D}(\boldsymbol{P})^{1}$. Show that $\boldsymbol{d}\{\boldsymbol{f}
,\boldsymbol{F}\}\in\operatorname{im}\boldsymbol{\Omega}^{1}$ as well and then apply
Proposition \ref{GauInv} again. Indeed,
\begin{align*}
\boldsymbol{d}\{\boldsymbol{f},\boldsymbol{F}\}  &  =\boldsymbol{d}%
\mathcal{L}_{\boldsymbol{X}}\boldsymbol{F}\\
&  =\mathcal{L}_{\boldsymbol{X}}\boldsymbol{dF}\\
&  =\mathcal{L}_{\boldsymbol{X}}i_{\boldsymbol{Z}}\boldsymbol{\omega}\\
&  =[\mathcal{L}_{\boldsymbol{X}},i_{\boldsymbol{Z}}]\boldsymbol{\omega
}+i_{\boldsymbol{Z}}\mathcal{L}_{\boldsymbol{X}}\boldsymbol{\omega}\\
&  =-i_{[\boldsymbol{X},\boldsymbol{Z}]}\boldsymbol{\omega+}i_{\boldsymbol{Z}%
}\boldsymbol{d}i_{\boldsymbol{X}}\boldsymbol{\omega}\\
&  =\boldsymbol{\Omega}^{1}(-[\boldsymbol{X},\boldsymbol{Z}%
])+i_{\boldsymbol{Z}}\boldsymbol{ddf}\\
&  =\boldsymbol{\Omega}^{1}([\boldsymbol{Z},\boldsymbol{X}]).
\end{align*}

\end{proof}

It is easy to prove that the action of $\boldsymbol{C}^{\infty}(\boldsymbol{P}%
)^{n-1}$ on gauge invariant elements in $\boldsymbol{C}^{\infty}%
(\boldsymbol{P})^{n}$ is indeed a Lie-algebra representation.

\begin{remark}
Notice that if the Euler-Lagrange equations are irreducible, then
$\boldsymbol{\Omega}$ is an isomorphism, $\ker\boldsymbol{\Omega}=0$ and every
element in $\boldsymbol{C}^{\infty}(\boldsymbol{P})^{\bullet}$ is trivially
gauge invariant. In this case $(\boldsymbol{C}^{\infty}(\boldsymbol{P}%
)^{n-1},\{{}\cdot{},{}\cdot{}\})$ acts on the whole $\boldsymbol{C}^{\infty
}(\boldsymbol{P})^{n}$.
\end{remark}

In \cite{bhs91} (see also \cite{fr05}) it has been shown that the bracket
described in full rigour in this section coincides with the Peierls bracket
\cite{p52}. In its turn the Peierls bracket is at the basis of a covariant
approach to quantization of field theories \cite{d03}. It is likely that the
mathematically rigorous picture presented here will help to better
understand, deal with and, possibly, generalize this complicated
\textquotedblleft functional\textquotedblright\ structure.

\subsection{Perspectives: Secondary Symplectic Reduction\label{Sec11}}

Most of the remarks in this section will be informal. From the physical point
of view, gauge invariant functions on $\boldsymbol{P}$ are the true
observables of the Lagrangian theory and, therefore, play a special role. We
shew in the last section that, basically, a Lie bracket is defined on gauge
invariant functions. We may go even further and ask:

\begin{enumerate}
\item are gauge invariant functions secondary functions on some secondary
manifold $\widetilde{\boldsymbol{P}}$?

\item if yes, is $\widetilde{\boldsymbol{P}}$ a symplectic reduction of the
secondary \textquotedblleft presymplectic manifold\textquotedblright%
\ $(\boldsymbol{P},\boldsymbol{\omega})$?
\end{enumerate}

In some more details, asking the last question amounts to wonder if there is
an embedding of algebras $\boldsymbol{\pi}^{\ast}:\boldsymbol{\Lambda
}(\widetilde{\boldsymbol{P}})^{\bullet}$ $\hookrightarrow\boldsymbol{\Lambda
}(\boldsymbol{P})^{\bullet}$ and a secondary two form $\widetilde
{\boldsymbol{\omega}}$ on $\widetilde{\boldsymbol{P}}$ such that 1)
$\ker\widetilde{\boldsymbol{\omega}}=0$ and 2) $\boldsymbol{\omega
}=\boldsymbol{\pi}^{\ast}(\widetilde{\boldsymbol{\omega}})$. Finding an answer
to the above questions would definitely establish the parallelism between
secondary calculus on the CPS and standard theory of constrained
(finite-dimensional) Hamiltonian systems. Moreover, it would fix the bases of
a mathematically rigorous, covariant, symplectic formalism for classical
Lagrangian field theories. Finally, it would represent a well founded starting
point for a covariant quantization of gauge systems \cite{ht92}.

A possible route through the answers to the above questions is described
below. First of all, there is a geometric counterpart of the degeneracy
distribution of $\boldsymbol{\omega}$. Let
\[
\xymatrix@C=30pt{0 \ar[r] & \varkappa |_{\mathscr{E}} \ar[r]^-{\ell_{\boldsymbol{E}(\mathscr{L})}} & \varkappa|_{\mathscr{E}}^\dag \ar[r]^-{\Delta_1} & P_2 \ar[r]^-{\Delta_2}  & \cdots}
\]
be a compatibility complex for $\ell_{\boldsymbol{E}(\mathscr{L})}$ and
\[
\xymatrix@C=30pt{0  & \varkappa |_{\mathscr{E}}^\dag \ar[l]  & \varkappa |_{\mathscr{E}} \ar[l]_-{\ell_{\boldsymbol{E}(\mathscr{L})}} & P_2^\dag \ar[l]_-{\Delta_{1}^\dag} & \cdots \ar[l]}
\]
its adjoint complex. There is an associated complex of $C^{\infty
}(\mathscr{E})$-modules:
\[
\xymatrix@C=30pt{0  & \overline{J}{}^\infty \varkappa |_{\mathscr{E}}^\dag \ar[l]  & \overline{J}{}^\infty\varkappa |_{\mathscr{E}} \ar[l]_-{h_{\boldsymbol{E}(\mathscr{L})}^\infty} & \overline{J}{}^\infty P_2^\dag \ar[l]_-{h_{1}^\infty} & \cdots \ar[l]},
\]
where we put $h_{\boldsymbol{E}(\mathscr{L})}^{\infty}:=h_{\ell
_{\boldsymbol{E}(\mathscr{L})}}^{\infty}$, $h_{1}^{\infty}:=h_{\Delta
_{1}^{\dag}}^{\infty}$ and so on. As discussed above, $\ker h_{\boldsymbol{E}%
(\mathscr{L})}^{\infty}\subset\overline{J}{}^{\infty}\varkappa|_{\mathscr{E}}$
identifies with $V\mathrm{D}(\mathscr{E})\subset V\mathrm{D}(J^{\infty
})|_{\mathscr{E}}$ via the isomorphism $\overline{J}{}^{\infty}\varkappa
|_{\mathscr{E}}\simeq V\mathrm{D}(J^{\infty})|_{\mathscr{E}}$ that sends
$\overline{j}{}_{\infty}\chi$ to $\rE_{\chi}$ (see Section \ref{HorCalc}). In
particular, $\ker h_{\boldsymbol{E}(\mathscr{L})}^{\infty}$ has got a natural
Lie algebra structure. Similarly, $\operatorname{im}h_{1}^{\infty}$
identifies with the module of sections of an involutive distribution
$\mathscr{G}$ on $\mathscr{E}$ made of vertical vector fields.

\begin{proposition}
\label{InvDist}$\operatorname{im}h_{1}^{\infty}\subset\ker h_{\boldsymbol{E}%
(\mathscr{L})}^{\infty}$ is a Lie-subalgebra.
\end{proposition}

\begin{proof}
Let $j_{1},j_{2}\in\overline{J}{}^{\infty}P_{2}^{\dag}$. Then $j_{1}=\sum
f_{1}\overline{j}{}_{\infty}\vartheta_{1}$ and $j_{2}=\sum f_{2}\overline{j}%
{}_{\infty}\vartheta_{2}$ for some $\ldots,f_{1},f_{2},\ldots\in C^{\infty
}(J^{\infty})$ and $\ldots,\vartheta_{1},\vartheta_{2},\ldots\in P_{2}^{\dag}$
(see Footnote \ref{page}, Section \ref{Sec9}, p.~\pageref{page}). Moreover,
$h_{1}^{\infty}(j_{1}),h_{1}^{\infty}(j_{2})$ correspond to vector fields
$X_{1}:=\sum f_{1}\rE_{\Delta_{1}^{\dag}\vartheta_{1}},X_{2}:=\sum
f_{2}\rE_{\Delta_{1}^{\dag}\vartheta_{2}}$, respectively, via the isomorphism
$\ker h_{\boldsymbol{E}(\mathscr{L})}^{\infty}\simeq V\mathrm{D}%
(\mathscr{E})$. Compute
\begin{align*}
\lbrack X_{1},X_{2}]  &  =\sum[f_{1}\rE_{\Delta_{1}^{\dag}\vartheta_{1}}%
,f_{2}\rE_{\Delta_{1}^{\dag}\vartheta_{2}}]\\
&  =\sum(f_{1}(\rE_{\Delta_{1}^{\dag}\vartheta_{1}}f_{2})\rE_{\Delta_{1}%
^{\dag}\vartheta_{2}}-f_{2}(\rE_{\Delta_{1}^{\dag}\vartheta_{2}}%
f_{1})\rE_{\Delta_{1}^{\dag}\vartheta_{1}}+f_{1}f_{2}\rE_{\{\Delta_{1}^{\dag
}\vartheta_{1},\Delta_{1}^{\dag}\vartheta_{2}\}}).
\end{align*}
Now, recall that $\mathfrak{g}=\operatorname{im}\Delta_{1}^{\dag}\subset
\ker\ell_{\boldsymbol{E}(\mathscr{L})}$ is a Lie subalgebra (see Corollary
\ref{CorsubLie}) so that $\{\Delta_{1}^{\dag}\vartheta_{1},\Delta_{1}^{\dag
}\vartheta_{2}\}=\Delta_{1}^{\dag}\vartheta$ for some $\vartheta\in
P_{2}^{\dag}$. Put
\[
j:=\sum f_{1}(\rE_{\Delta_{1}^{\dag}\vartheta_{1}}f_{2})\overline{j}_{\infty
}\vartheta_{2}-f_{2}(\rE_{\Delta_{1}^{\dag}\vartheta_{2}}f_{1})\overline
{j}_{\infty}\vartheta_{1}+f_{1}f_{2}\overline{j}_{\infty}\vartheta\in
\overline{J}{}^{\infty}P_{2}^{\dag}.
\]
Then $h_{1}^{\infty}(j)\in\operatorname{im}h_{1}^{\infty}\subset\ker
h_{\boldsymbol{E}(\mathscr{L})}^{\infty}$ corresponds to $[X_{1},X_{2}]\in
V\mathrm{D}(\mathscr{E})$ via the isomorphism $\ker h_{\boldsymbol{E}%
(\mathscr{L})}^{\infty}\simeq V\mathrm{D}(\mathscr{E})$.
\end{proof}

In the following we will understand isomorphism $\ker h_{\boldsymbol{E}%
(\mathscr{L})}^{\infty}\simeq V\mathrm{D}(\mathscr{E})$. In view of
Proposition \ref{InvDist}, $\mathscr{G}$ is an involutive distribution on
$\mathscr{E}$. Namely, $\mathscr{G}$ is the (involutive) distribution
generated by evolutionary derivatives with generating sections in
$\mathfrak{g}$ (such kinds of distributions have been recently considered in
\cite{kv08}).

Notice that the horizontal Spencer differential $\overline{S}:V\mathrm{D}%
(\mathscr{E})\otimes\overline{\Lambda}(\mathscr{E})\longrightarrow$
$V\mathrm{D}(\mathscr{E})\otimes\overline{\Lambda}(\mathscr{E})$
\textquotedblleft restricts\textquotedblright\ to $\operatorname{im}%
h_{1}^{\infty}\otimes\overline{\Lambda}(\mathscr{E})$. Denote by
$\overline{s}:\operatorname{im}h_{1}^{\infty}\otimes\overline{\Lambda
}(\mathscr{E})\longrightarrow\operatorname{im}h_{1}^{\infty}\otimes
\overline{\Lambda}(\mathscr{E})$ the restricted differential. Clearly,
$\mathfrak{g}\subset\mathfrak{g}_{1}:=H^{0}(\operatorname{im}h_{1}^{\infty
}\otimes\overline{\Lambda}(\mathscr{E}),\overline{s})=\mathbf{D}%
(\boldsymbol{P})^{0}\cap\operatorname{im}h_{1}^{\infty}$. We now describe the
quotient $\mathfrak{g}_{1}/\mathfrak{g}$. Let $\square:\varkappa
|_{\mathscr{E}}\longrightarrow Q_{2}$ be a compatibility operator for
$\Delta_{1}^{\dag}:P_{2}^{\dag}\longrightarrow$ $\varkappa|_{\mathscr{E}}$,
and put $k:=h_{\square}^{\infty}:\overline{J}{}^{\infty}\varkappa
|_{\mathscr{E}}\longrightarrow\overline{J}{}^{\infty}Q_{2}$. Then
$\operatorname{im}h_{1}^{\infty}=\ker k$, so that
\[
\mathfrak{g}_{1}=H^{0}(\operatorname{im}h_{1}^{\infty}\otimes\overline
{\Lambda}(\mathscr{E}),\overline{s})=H^{0}(\ker k\otimes\overline{\Lambda
}(\mathscr{E}),\overline{s})=\ker\square.
\]
We conclude that $\mathfrak{g}_{1}/\mathfrak{g}=\ker\square/\operatorname{im}%
\Delta_{1}^{\dag}\simeq H^{1}(\operatorname{im}h_{1}^{\infty}\otimes
\overline{\Lambda}(\mathscr{E}),\overline{s})$ (see Theorem \ref{SpencTh}) and
there is an exact sequence (of vector spaces)
\[
\xymatrix{0\ar[r] & \mathfrak{g}\ \ar@{^(->}[r] & \mathfrak{g}_1%
\ar@{->>}[r] & H^{1}(\operatorname{im}h_{1}^{\infty}\otimes\overline{\Lambda
}(\mathscr{E}),\overline{s})\ar[r] & 0}.
\]
Thus $H^{1}(\operatorname{im}h_{1}^{\infty}\otimes\overline{\Lambda
}(\mathscr{E}),\overline{s})$ is the obstruction to $\mathfrak{g}$ being
isomorphic to $0$-cohomology of the complex $(\operatorname{im}h_{1}^{\infty
}\otimes\overline{\Lambda}(\mathscr{E}),\overline{s})$, that is, in a sense,
the obstruction to \textquotedblleft the algebraic description and the
geometric description of gauge symmetries coinciding\textquotedblright.

Despite the possible existence of such an obstruction, define the new
distribution on $\mathscr{E}$, $\widetilde{\mathscr{C}}%
:=\mathscr{E}+\mathscr{G}$. $\widetilde{\mathscr{C}}$ is, generally,
infinite-dimensional. Moreover, it is an involutive distribution. Roughly
speaking, integral submanifolds of $\widetilde{\mathscr{C}}$ identify with
\textquotedblleft gauge equivalence classes\textquotedblright\ of solutions of
the Euler-Lagrange equations. Therefore, it is natural to put $\widetilde
{\boldsymbol{P}}:=\{$maximal integral submanifolds of $\widetilde
{\mathscr{C}}\}$ and interpret $\widetilde{\boldsymbol{P}}$ as the space of
\textquotedblleft physical states\textquotedblright\ of fields of the
Lagrangian theory $(\pi,\boldsymbol{S})$.

A secondary calculus may be introduced on $\widetilde{\boldsymbol{P}}$,
basically via the $\widetilde{\mathscr{C}}$-spectral sequence $\widetilde
{\mathscr{C}}E(\mathscr{E})$, so that elements in $\widetilde{\mathscr{C}}%
E_{1}(\mathscr{E})=:\boldsymbol{\Lambda}(\widetilde{\boldsymbol{P}})^{\bullet
}$ are interpreted as (secondary) differential forms on $\widetilde
{\boldsymbol{P}}$. The inclusion $\mathscr{C}\subset\widetilde{\mathscr{C}}$
induces a morphism $\widetilde{\mathscr{C}}E(\mathscr{E})\longrightarrow
\mathscr{C}E(\mathscr{E})$ of spectral sequences whose $1$-st term we denote
by $\boldsymbol{\pi}^{\ast}:\boldsymbol{\Lambda}(\widetilde{\boldsymbol{P}%
})^{\bullet}\longrightarrow\boldsymbol{\Lambda}(\boldsymbol{P})^{\bullet}$.

Now, we'd like to interpret $\widetilde{\boldsymbol{P}}$ as a
\textquotedblleft(symplectically) reduced CPS\textquotedblright. In order to
be able to do this in a consistent and physically meaningful way at least the
following two conditions should be fulfilled:

\begin{enumerate}
\item the image of $\boldsymbol{C}^{\infty}(\widetilde{\boldsymbol{P}%
})^{\bullet}:=\widetilde{\mathscr{C}}E_{1}^{0,\bullet}(\mathscr{E})$ under
$\boldsymbol{\pi}^{\ast}$ should be made of gauge invariant (secondary)
functions on $\boldsymbol{P}$,

\item a secondary $2$-form $\widetilde{\boldsymbol{\omega}}$ on
$\widetilde{\boldsymbol{P}}$ should exist so that $\ker\widetilde
{\boldsymbol{\omega}}=0$ and $\boldsymbol{\pi}^{\ast}(\widetilde
{\boldsymbol{\omega}})=\boldsymbol{\omega}$.
\end{enumerate}

If this was the case then, in the author's opinion, $(\widetilde
{\boldsymbol{P}},\widetilde{\boldsymbol{\omega}})$ could be \textquotedblleft
safely\textquotedblright\ referred to as the \textquotedblleft symplectic
reduction of $(\boldsymbol{P},\boldsymbol{\omega})$\textquotedblright\ from
the mathematical point of view, and as the \textquotedblleft reduced
CPS\textquotedblright\ \cite{ht92,lw90,r04} from the physical point of view.

As suggested by the example in this section and by preliminary work by the
author, typical homological algebra (and, possibly, homological perturbation
theory) techniques seem to be necessary to investigate further in this
direction and complete the above sketched program.

\section*{Conclusions}

We proposed a fully rigorous approach to the geometry of the covariant phase
space $\boldsymbol{P}$, and the canonical, closed $2$-form
$\boldsymbol{\omega}$ on it, in the framework of secondary calculus. In
particular, we described the kernel of $\boldsymbol{\omega}$ in terms of the
compatibility operator for the linearized Euler-Lagrange equations thus
revealing the precise relation between gauge symmetries and constraints in
field theory \cite{lw90}. We also described gauge invariant (secondary)
functions on $\boldsymbol{P}$ and their Lie algebra structure. It is likely
that such a Lie algebra is at the basis of a covariant canonical quantization
of the theory \cite{d03}. A step forward in this direction would be to
rigorously perform a symplectic reduction of $(\boldsymbol{P}%
,\boldsymbol{\omega})$. The preliminary analysis presented in Section
\ref{Sec10} suggests that this is possible, and should be done, within
secondary calculus (or a slight generalization of it) and, in any case, by
means of cohomological techniques.

We stress that, in this paper, we basically worked \textquotedblleft on
shell\textquotedblright. The relationship with \textquotedblleft off
shell\textquotedblright\ methods (Koszul-Tate resolution and BRST complex
\cite{bbh95,bbh00,ht92} - see also \cite{v02}) should be carefully analyzed.

\end{document}